\documentclass[12pt,reqno]{amsart}
\usepackage{latexsym}
\usepackage{amsmath}
\usepackage{amssymb}
\usepackage{amsfonts}
\usepackage{verbatim}
\usepackage[latin1]{inputenc}
\usepackage{color}

\newtheorem{Th}{Theorem}[section]
\newtheorem{prop}[Th]{Proposition}
\newtheorem{Df}[Th]{Definition}

\newtheorem{Ex}[Th]{Example}
\newtheorem{Def}[Th]{Definition}
\newtheorem{lem}[Th]{Lemma}
\newtheorem{rk}[Th]{Remark}

\let \ssection=\section
\renewcommand{\section}{\setcounter{equation}{0}\ssection}

\def\^#1{\if#1i{\accent"5E\i}\else{\accent"5E #1}\fi}
\def\"#1{\if#1i{\accent"7F\i}\else{\accent"7F #1}\fi}

\def\^#1{\if#1i{\accent"5E\i}\else{\accent"5E #1}\fi}
\def\"#1{\if#1i{\accent"7F\i}\else{\accent"7F #1}\fi}
\allowdisplaybreaks

\begin{document}
\title[Discretization schemes]
{On discretization schemes for stochastic evolution equations}

\author[I. Gy\"ongy]
{Istv\'an Gy\"ongy}
\thanks{This paper was written while the first named author was visiting
the University of Paris X. The research of this author is partially supported by
EU Network HARP}
\address{School of Mathematics,
University of Edinburgh,
King's  Buildings,
Edinburgh, EH9 3JZ, United Kingdom}
\email{gyongy@maths.ed.ac.uk}

\author[A. Millet]{Annie Millet}
\thanks{The research of the second named author is partially supported
by the research project BMF2003-01345}
\address
{Laboratoire de Probabilit\'es et Mod\`eles Al\'eatoires
(CNRS UMR 7599), Universit\'es Paris~6-Paris~7, Boite Courrier 188, 4 place Jussieu, 75252 Paris Cedex 05,
{\it  and } SAMOS-MATISSE,
Universit\'e Paris 1, 90 Rue de Tolbiac, 75634 Paris Cedex 13}
\email{amil@ccr.jussieu.fr}



\begin{abstract}
Stochastic evolutional equations with monotone
operators are considered in Banach spaces. Explicit and implicit
numerical schemes are presented. The convergence of the
approximations to the solution of the equations is proved.
\end{abstract}

\maketitle

\section{Introduction}
                                            \label{intro}
Let $V\hookrightarrow  H\hookrightarrow   V^*$ be a
{\it normal triple} of spaces with dense and
continuous embeddings, where $V$
is a reflexive Banach space,
$H$ is a Hilbert space, identified with its dual by means 
of the inner product in $H$, and $V^*$ is the dual of $V$.  Let
$W=(W_t)_{t\geq 0}$ be an $r$-dimensional Brownian motion carried by a stochastic basis 
$(\Omega, \mathcal F, ({\mathcal F}_t)_{t\geq 0}, P)$. In this paper, we study the approximation
 of the solution to the evolution equation
\begin{equation}\label{u}
 u_t=u_0 + \int_0^t A_s(u_s)\, ds +  \sum_{j=1}^r \int_0^t
 B_s^j(u_s)\, dW^j_s\, ,
\end{equation}
where $u_0$ is a $H$-valued ${\mathcal F}_0$-measurable random variable,
$A$ and $B$ are  (non-linear) adapted  operators  defined on $[0,+\infty[\times V\times \Omega$
with values in $V^*$ and $H^r$ respectively.  
  
The conditions imposed on $A_s$ are satisfied
by the following classical
example: $V=W^{1,p}_0(D),\;  H=L^2(D) \; V^*=W^{-1,q}(D)$ and 
\[ A_s(u)={\displaystyle \sum_{i=1}^d }
\frac{\partial }{\partial x_i}\,\left(  \left| \frac{\partial u}{\partial x_i}\right|^{p-2}\,
 \frac{\partial u}{\partial x_i}\right)\, , \]
where $D$ is a bounded domain of
${\mathbb R}^d$, $p\in ]2,+\infty[$ and $q$ are conjugate
exponents.  In \cite{L} the monotonicity method is used in the
deterministic case to prove that if $u_0\in H$ and $B=0$, equation
(\ref{u}) has a unique solution in $L^p_V(]0,T])$ such that
$u_t=0$ on  $]0,T]\times \partial D$.
 Using the monotonicity method,
the existence and uniqueness of a solution  $u$  to (\ref{u}) is
proved in \cite{P} and \cite{KR}.  This result can be fruitfully
applied also to linear stochastic PDEs, in particular to the
equations of nonlinear filtering theory (see \cite{Pa}, \cite{Par}
and \cite{R}).  The existence and uniqueness theorem from \cite{KR}
is extended  in \cite{G} to equation (\ref{u}) with martingales
and martingale measures in place of $W$. Inspired by \cite{K}, the
method of monotonicity is interpreted in \cite{GM} as a
minimization method for some convex functionals.

In the present paper we introduce an implicit time discretization
$u^m$,  space-time explicit and implicit discretization schemes
$u^m_n$ and $u^{n,m}$ of $u$ defined in terms of a constant time
mesh ${\delta_m}  = \frac{T}{m}$ and of a sequence of finite
dimensional subspaces  $V_n$ of $V$.   One particular case of such
spaces is that used  in the Galerkin method or in the
piecewise linear finite elements methods.  To define space-time
discretizations of $u$, we denote by $\Pi_n : V^* \rightarrow V_n$
a $V_n$-valued projection.
\smallskip

For $0\leq i\leq m$, set $t_i=\frac{i\, T}{m}$.
The  explicit $V_n$-valued
space-time discretization of $u$
is defined for an initial condition $u_0\in H$  by
$u^n_m(t_0)= u^n_m(t_1)=\Pi_nu_0$  and for $1\leq i<m$,
\begin{equation}                                                    \label{spacetimee}
u^n_m(t_{i+1})= u^n_m(t_i) +
{\delta_m}  \, \Pi_n\tilde{A}^m_{t_i}(u^n_m(t_i)) 
  +\sum_{j=1}^r \Pi_n\tilde{B}^{m,j}_{t_i}(u^n_m(t_i))\,
\big(W^j_{t_{i+1}}- W^j_{t_i}\big)\, ,
\end{equation}
where for $x\in V$, $\tilde{A}^m_{t_i}(x)\in V^*$
and
$(\tilde{B}^{m,j}_{t_i}(x)\, ,\, 1\leq j\leq r)\in H^r$
denote the averages
of the processes $A_.(x)$ and $B_.(x)$
over the time interval $[t_{i-1},t_i]$.
\smallskip

The $V$-valued implicit time discretization of
$u$ is defined  for an initial condition
$u_0\in H$ by $u^m(t_0)=0$,
$u^m(t_1)=u_0+{\delta_m}  \, A^m_{t_1}(u^m(t_{1}))$,
and for $1\leq i<m$,
\begin{equation}\label{time1}
u^m(t_{i+1})= u^m(t_i) + {\delta_m}  \, A^m_{t_i}(u^m(t_{i+1}))
+ \sum_{j=1}^r \tilde{B}^{m,j}_{t_i}(u^m(t_i))\,
\big(W^j_{t_{i+1}}- W^j_{t_i}\big)\, ,
\end{equation}
where for $x\in V$, $A^m_{t_i}(x)\in V^*$ denotes  the average of
the process $A_.(x)$ over the time interval $[t_i,t_{i+1}]$ and as
above $(\tilde{B}^{m,j}_{t_i}(x)\, ,\, 1\leq j\leq r )\in H^r$
denotes the average of $B_.(x)$ over the time interval $[t_{i-1},
t_i]$.
\smallskip

\noindent Finally, the implicit $V_n$-valued
space-time discretization of $u$ is defined for
$u_0\in H$  by $u^{m,m}(t_0)=0$,
$u^{n,m}(t_1)=\Pi_nu_0+{\delta_m}
\, \Pi_n A^m_{t_1}(u^{n,m}(t_{1}))$,
and for  $1\leq i<m$,
\begin{equation}                                              \label{spacetimei}
u^{n,m}(t_{i+1})= u^{n,m}(t_i) + {\delta_m}  \,
\Pi_n A^m_{t_i}(u^{n,m}(t_{i+1}))  
 + \sum_{j=1}^r \Pi_n \tilde{B}^{m,j}_{t_i}(u^{n,m}(t_i))\,
\big(W^j_{t_{i+1}}- W^j_{t_i}\big)\, ,
\end{equation}
where $A^m_{t_i}$  and $\tilde{B}^{m,j}_{t_i}$
have been defined above.
\smallskip 

The processes $v$ equal to $u^m$, $u^n_m$  or
$u^{n,m}$ are defined between $t_i$ and $t_{i+1}$
as stepwise
constant adapted stochastic processes, i.e.,
$v(t):=v(t_i)$ for $t\in ]t_i, t_{i+1}[$.
We prove that for $m$ large enough, (\ref{time1})
(resp. (\ref{spacetimei}))
has a unique solution $u^m$  (resp. $u^{n,m}$),
which  converges  weakly to $u$
in a weighted space of $p$-integrable processes, and that
the approximations at terminal time $T$
converge strongly to $u(T)$  in
 $L^2_H(\Omega)$ as $m\rightarrow +\infty$
(resp. $n$ and $m$ go to infinity).  
As one expects, the  convergence of the
explicit approximation $u^n_m$ to $u$ in these
spaces requires some condition
relating the time mesh $T/m$ and the spaces $V_n$.
The existence of the solution to
(\ref{time1}) or (\ref{spacetimei}), as well
as that of a limit for some subsequence $u^{m_k}$,
$u^{n_k}_{m_k}$ or   $u^{n_k,m_k}$ is proved using
apriori estimates,
which are  based on the coercivity, monotonicity
and growth assumptions made on the operators
$A_s$ and $B_s$.  The identification of $u$ as
the limit is obtained by means of
the minimization property
of $u$.  Note that the conditions imposed
on the operators $A_s$ and $B_s$ involve
constants which may depend on time.
This  allows the operators to approach degeneracy. 
However, this lack of uniform non-degeneracy has to be
balanced by a suitable growth condition
which depends on time as well.  
Thus, as a by-product of the identification of the weak
limit  of the explicit and
implicit space-time discretization schemes,
we obtain the existence of a solution
to (\ref{u}) under slightly more general conditions than 
those used in \cite{Pa}, \cite{KR} or \cite{GM}.
\smallskip 

Section 2 states the conditions imposed on the operators $A$ and
$B$, the spaces $V_n$ and the maps $\Pi_n$, gives examples
satisfying these conditions, describes precisely the explicit and
implicit  schemes, and states the corresponding convergence
results.  The third section provides the proofs of the main
theorems  and an appendix gathers some technical tools.
\smallskip

As usual we denote by $C$ a constant
which may change from line to line.  
All the processes considered will be adapted
with respect to the filtration
$({\mathcal F}_t\, ,\, t\geq 0)$.
\section{Description of the results}
We first state the precise assumptions made on the  operators.
Let $V$ be a separable reflexive Banach space,
 embedded continuously and densely
into a Hilbert space $H$, which is identified with its dual,
$H^*$ by means 
of the inner product $(\cdot,\cdot)$ in $H$.
Then the adjoint embedding
$H\hookrightarrow V^*$ of $H^*\equiv H$ into $V^{*}$, the dual of $V$,
is also dense and continuous.
Let  $\langle v,x\rangle=\langle x,v\rangle$ denote the duality product
for $v\in V$ and $x\in V^*$.
Observe that $\langle v,h\rangle =(v,h)$
for $h\in H$ and $v\in V$.  Let $(\Omega, {\mathcal  F}, ({\mathcal F}_t)_{t\geq 0}, P)$ be a
stochastic basis, satisfying the usual conditions and carrying an
$r$-dimensional Wiener martingale $W=(W_t)_{t\geq0}$ with respect
to $({\mathcal F_t})_{t\geq 0}$.
 \smallskip

Fix $T>0$, $p\in [2,+\infty[$ and let $q=\frac{p}{p-1}$  be  the
conjugate exponent of $p$. Let $L^1$ (resp. $L^2$) denote the
space of integrable (resp. square integrable)
 real functions  over $[0,T]$.  
Let
\[ A:[0,T] \times V\times \Omega \rightarrow V^*\,
, \quad B:[0,T]\times V\times \Omega \rightarrow H^r\]
be such that for every $v,w\in V$ and
$1\leq j\leq r$, $\langle w,A_s(v)\rangle $ and
$(B_s^j(v),w)$ are adapted processes and
the following conditions hold:
\medskip

{\bf (C1)} The pair $(A,B)$ satisfies the
{\it monotonicity condition},   i.e.,
almost surely for all  $t\in [0,T]$,
$x$ and $y$ in $V$:
\begin{equation}\label{monotone}
2\, \langle x-y, A_t(x)-A_t(y)\rangle
+ \sum_{j=1}^r |B_t^j(x)-B_t^j(y)|^2_H \leq 0\, .
\end{equation}
\indent {\bf (C2)} The pair $(A,B)$ satisfies the {\it coercivity
condition}  i.e., there exist  non-negative 
integrable  functions  $K_1,\bar
K_1$ and $\lambda:]0,T]\rightarrow  ]0,+\infty[$ 
 such that almost surely
\begin{equation}        \label{coervive}
2\, \langle x,A_t(x)\rangle  + \sum_{j=1}^r |B^j_t(x)|^2_H
+ \lambda(t)\, |x|^p_V \leq K_1(t)|x|^2_H+\bar K_1(t)
\end{equation}
for all $t \in  ]0,T]$ and $x\in V$.

{\bf (C3)} The operator $A$ is {\it hemicontinuous}  i.e.,
almost surely
\begin{equation}                               \label{hemi}
\lim_{\varepsilon \rightarrow 0}
\langle A_t(x+\varepsilon y), z\rangle
= \langle A_t(x),z\rangle \, .
\end{equation}
for all  $t\in [0,T]$, $x,y,z$ in $V$.

{\bf (C4)}
({\it  Growth condition})
There exist a non-negative function $K_2\in L^1$
and a constant $\alpha\geq 1$ such that almost surely
\begin{equation}                        \label{growthA}
|A_t(x)|^q_{V^*}
\leq \alpha \lambda^q(t)\,|x|^{p}_V
+ \lambda^{q-1}(t){K}_2(t)
\end{equation}
for all  $t\in ]0,T]$ and $x\in V$.

We also impose some integrability of the initial condition $u_0$:
\medskip

{\bf (C5)} $u_0\, :\, \Omega\rightarrow H$ is
${\mathcal F}_0$-measurable and such that
$E(|u_0|_H^2)<+\infty$.
\medskip

\begin{rk}                                   \label{15.29.10}
From {\bf (C2)} and {\bf (C4)} it
is easy to get that almost surely
\begin{equation}                                \label{10.30.11}
\sum_{j=1}^r|B^j_t(x)|^2_H\leq (2\, \alpha +1) \,  \lambda(t)|x|_{V}^p
+K_1(t)|x|^2_H+ K_3(t)
\end{equation}
for all $t\in ]0,T]$ and $x\in V$, where
 $K_3(t)=\bar{K}_1(t)+\frac{2}{q}\, K_2(t) \in L^1$.
\end{rk}
\begin{proof}
For every $t\in ]0,T]$ and $x\in V$,
\begin{eqnarray*}  
 |\langle x\, ,\, A_t(x)\rangle|& \leq  & |x|_V\, |A_t(x)|_{V^*} \leq
\alpha^{\frac{1}{q}}\, \lambda(t)\, |x|_V^{1+\frac{p}{q}} + \lambda(t)^{\frac{q-1}{q}}\, 
|x|_V\, K_2(t)^{\frac{1}{q}}\\
 &\leq& \alpha^{\frac{1}{q}}\, \lambda(t) |x|_V^p 
+ \frac{1}{p}\, \lambda(t)\, |x|_V^p + \frac{1}{q}\, K_2(t)\, .
\end{eqnarray*}
Thus, (\ref{coervive}) and (\ref{growthA}) yield (\ref{10.30.11}).
\end{proof}
Note that, unlike in \cite{GM}, \cite{KR} and \cite{Pa}, 
the coercivity constant $\lambda(t)$ can vary with $t$
(for example, one can suppose that
$\lambda(t) = \lambda\, t$ for some constant $\lambda >0$), which means that the operators
can be more and more degenerate as $t\rightarrow 0$.
 However, this bad behavior has to be balanced by some
more and more stringent growth conditions.

We remark that the monotonicity condition (C1) can be weakened as follows:
\smallskip

{\bf (C1bis)}  There exists a non negative function $K\in L^1_+$ such
that almost every $(t,\omega)\in[0,T]\times\Omega$ and every
$x,y\in V$
$$
 2\, \langle x-y, A_t(x)-A_t(y)\rangle
+ \sum_{j=1}^r |B_t^j(x)-B_t^j(y)|^2_H \leq K(t)\, |x-y|_H^2 \, .
$$

\noindent Indeed, if $u$ is  a  solution to (\ref{u})  and
$\gamma_t:=\exp\Big(\frac{1}{2}\int_0^tK(s)\,ds\Big)$, then
$v_t=\gamma^{-1}_t u_t$ is a solution of the equation
\[ v_t=u_0+\int_0^t\, \bar{A}_s(v_s)\, ds 
+\sum_{j=1}^r \int_0^t \bar{B}^j_s(v_s)\, dW^j_s\, ,\]
where for every $t\in[0,T]$ and $x\in V$:
\[
\bar{A}_t(x):=\gamma^{-1}_t\,
A_t\big(\gamma_t\, x\big) - \frac{1}{2}{K(t)}\, x,
\; \mbox{ \rm and } \;
\bar{B}_t(x):=\gamma^{-1}_t\
B_t\big(\gamma_t\,  x \big)\, . 
\]
If $(A,B)$ satisfies  (C1bis) then it is easy to see that
$(\bar{A},\bar{B})$ satisfies (C1).  Clearly, if $A$ is
hemicontinuous, then $\bar A$ is also hemicontinuous.  If $(A,B)$
satisfies the coercivity condition  (C2), then $(\bar A,\bar B)$
also satisfies (C2).  If $A$ satisfies the growth condition (C4)
then it is an easy exercise to check that $\bar{A}$ also satisfies
(C4), provided $p\geq 2$ and $K(t)\leq C\lambda(t)$ for all $t$
with some constant $C$.

\begin{Ex} A large class of linear 
and semi-linear stochastic 
partial differential equations of parabolic type satisfies the 
above conditions.  Below we present 
 a class of examples  
of nonlinear  equations.  Let $D$ be a bounded domain of 
${\mathbb R}^d$, 
$p\in [2,+\infty[$,
$V=W_0^{1,p}(D)$, $H=L^2(D)$, $V^*=W^{-1,q}(D)$. 
Let  the operators 
$A_t$, $B^j_t$ be defined by 
$$
 A_t(u,\omega)
:= \sum_{i=1}^d 
\frac{\partial}{\partial x_i}
f_i(t,x,\nabla u(x),\omega),
$$
$$
 B^k_t(u,\omega):=g^k(t,x,\nabla u(x),\omega)
+h^k(t,x,u(x),\omega), \quad k=1,2,...,r 
$$
for 
$u\in V$, $t\in[0,T]$ 
and $\omega\in\Omega$,  
where 
$\nabla u$ denotes the gradient of $u$, i.e., 
$\nabla u= (\frac{\partial u}{\partial
x_1},\frac{\partial u}{\partial x_2},
 ..., \frac{\partial u}{\partial x_d})  $, 
and
$f_i=f_i(t,x,z,\omega)$, 
$g^j=g^j(t,x,z,\omega)$, 
$h^j=h^j(t,x,s,\omega)$ 
are some real valued functions of 
$t\in[0,\infty[$, $x,z\in{\mathbb R}^d$ 
and $s\in\mathbb R$,  
such that the following conditions are 
satisfied:

(i)  The functions $f_i$, $g^j$ and $h^j$
are Borel measurable in $t,x,z,s$ for
each fixed $\omega$, and 
are ${\mathcal F}_t$-adapted stochastic processes 
for each fixed $t,x,z,s$.

(ii) The functions $f_i$ and $g^j$ are 
differentiable in $z=(z_1,z_2,...,z_d)$, and there 
exists a constant $\varepsilon>0$, such that for 
almost every $\omega\in\Omega$ and all $t,x,z$ the 
matrix 
$$
\big(S_{ij}\big)
:=\big(2f_{iz_j}
-(1+\varepsilon)\sum_{k=1}^rg^k_{z_i}g^k_{z_j}\big)
$$ 
is positive semidefinite, where 
$f_{iz_j}:=\frac{\partial}{\partial z_j}f_i$, 
$g^k_{z_j}:=\frac{\partial}{\partial z_j}g^k$.

(iii) There exists a function $K:[0,T]\to[0,\infty[$, 
$K\in L^1$, such that 
$$
\sum_{k=1}^r|h^k(t,x,u)-h^k(t,x,v)|^2
\leq K(t)|u-v|^2,
$$
$$
\sum_{k=1}^r\int_{{\mathbb R}^d}
|h^k(t,x,0)|^2\,dx
\leq K(t)
$$
for almost every $\omega\in\Omega$ and 
all $t\in[0,T], x\in{\mathbb R}^d$,  
$u,v\in\mathbb R$.

(iv) There exist a constant 
$\varepsilon>0$ and a function 
$\lambda:]0,T]\to]0,\infty[$, 
$\lambda\in L^1$, such that 
almost surely
$$
2\sum_{i=1}^dz_if_i(t,x,z)
-(1+\varepsilon)\sum_{k=1}^r|g^k(t,x,z)|^2
\geq \lambda(t)|z|^p, 
$$
$$
\sum_{i=1}^d|f_i(t,x,z)|
\leq \alpha\lambda(t)|z|^{p-1}
+\lambda^{\frac{1}{p}}(t)K_1^{\frac{1}{q}}(t,x)
$$
for all $t\in ]0,T]$, $x,z\in{\mathbb R}^d$, 
where $\alpha>0$ is a constant and 
$K_1:[0,T]\times{\mathbb R}^d\to[0,\infty[$ is 
a function such that for every $t\in ]0,T]$,
$\int_{{\mathbb R}^d}    K_1(t,x)\,dx<\infty$  and 
$\int_0^T\int_{{\mathbb R}^d}K_1(t,x)\,dx\,dt<\infty$. 

It is an easy exercise to verify that under these conditions 
$A$ and $(B^j)$ 
satisfy conditions (C2)-(C4) and (C1bis). A simple example 
of nonlinear functions $f_i$, $g^k$ and $h^k$, satisfying 
the above conditions (i)-(iv), is for $p\in ]2,+\infty[$
\begin{eqnarray*}
f_i(t,x,z,\omega)&:= & a_i(t,x,\omega)|z_i|^{p-2}z_i,\\
g^k(t,x,z,\omega)&:=&2p^{-1}\sum_{i=1}^d\, 
b_i^k(t,x,\omega)|z_i|^{\frac{p}{2}},\\
h^k(t,x,u,\omega)&:=&c^k(t,x,\omega)|u|+d^k(t,x,\omega)
\end{eqnarray*}
for $t\in[0,T]$, $x,z=(z_1,...,z_d)\in{\mathbb R}^d$, 
$u\in\mathbb R$, $\omega\in\Omega$, where $a_i$, $b^k_i$, 
$c^k$ and $d^k$ are real valued  functions such that the following 
conditions hold:

(1) The functions $a_i$, $b^k_i$, 
$c^k$ and $d^k$ are Borel functions of $t,x$ 
for
each fixed $\omega$, and 
are ${\mathcal F}_t$-adapted stochastic processes 
for each fixed $x$.

(2) There exist constants  
$\varepsilon>0$, $\alpha>0$ and a function 
$\lambda:]0,T]\to]0,\infty[$, 
$\lambda\in L^1$, such that 
almost surely
$$
\Big(2(p-1)\, a_i(t,x) \, \delta_{ij}  \, 
-(1+\varepsilon)\sum_{k=1}^r  (b^k_ib^k_j)(t,x) , 1 \leq i,j \leq d \Big)
\geq\lambda(t)I
$$
$$
\sum_{i=1}^d a_i(t,x) \leq \alpha\lambda(t)
$$
for all $t\in ]0,T]$ and $x\in{\mathbb R}^d$, where 
$I$ is the identity matrix, and $\delta_{ij}=1$ for $i=j$ and  
$\delta_{ij}=0$ otherwise. 

(3) There exist functions $ K : [0,T]\to [0,\infty[$ and 
$L:[0,T]\times{\mathbb R}^d\to[0,\infty[$ such that almost surely 
$$
\sum_{k=1}^r|c^k(t,x,\omega)|^2\leq K(t), 
\quad \sum_{k=1}^r|d^k(t,x,\omega)|^2\leq L(t,x)
$$
for all $t,x$, and 
$$
\int_0^T K(t)\, dt <\infty \, , \quad 
\int_0^T\int_{{\mathbb R}^d}L(t,x)\,dx\,dt<\infty.
$$

We remark that though for $p=2$ the function 
$g^k(t,x,z,\omega):= 
 \sum_{i=1}^d
b_i^k(t,x,\omega)|z_i|$ is not differentiable 
at points $z$ such that $z_i=0$ for some $i$, it is easy to see that the 
corresponding operators 
$A$, $B^k$ still satisfy conditions (C2)-(C4) and 
(C1bis) also in this case. 
\end{Ex}

Note that the conditions (C2)-(C4) slightly extend those used in
\cite{Pa}, \cite{KR} or \cite{GM}, where the function 
$\lambda$ is supposed to be constant.

\begin{Def}\label{sol}
An adapted continuous $H$-valued process $u$ is a solution to (\ref{u}) if

(i) $E\int_0^T  |u_t|_V^p \, \lambda(t)\, dt <\infty$.

(ii)  For every $t\in [0,T]$ and $z \in V$
\begin{equation}\label{solu}
\langle u_t, z\rangle  = \langle u_0,z\rangle
+ \int_0^t \langle A_s(u_s),z
\rangle \, ds
+ \sum_{j=1}^r \int_0^t (B_s^j(u_s),z)\,
dW^j_s\quad a.s.
\end{equation}
\end{Def}


Notice that under condition  (C4)  and (\ref{10.30.11}), i.e., for
example under conditions (C2) and (C4), it is easy to see that an
adapted continuous $H$-valued $u$ is a solution to (\ref{u}) as
soon as (\ref{solu}) is satisfied for all $z$ in a dense subset of
$V$.  The following theorem extends the existence and uniqueness
theorem proved in \cite{Pa} and \cite{KR}.

\begin{Th}                                                    \label{existu} 
Let conditions  (C1)-(C5) hold.  
Then equation (\ref{u}) has a
unique solution $u$.
\end{Th}

\begin{rk}\label{uniqueness}
The uniqueness of the solution to equation (\ref{u})
follows easily from conditions (C1) and (C4).   Moreover, 
if $u$ is a solution of equation (\ref{u}), 
then conditions (C2) and (C5) imply 
\begin{equation}                                              \label{becsles}
\sup_{t\in[0,T]}E|u_t|^2_H<\infty.
\end{equation}
\end{rk}

{\it Proof of Remark \ref{uniqueness}}.
Let $u^{(1)}$ and $u^{(2)}$ be  solutions to (\ref{u}).  Then for
$\delta_t:=u^{(1)}_t-u^{(2)}_t$ we have
\begin{equation}                                                   \label{ito}
\delta_t=\int_0^t z^{\ast}_s\,dY_s+h_t 
, \quad dY_t\times dP -a.e.,
\end{equation}
where

\begin{eqnarray*}
z^{\ast}_t:&=&\lambda^{-1}(t)
\big[A_t(u^{(1)}_t)-A_t(u^{(2)}_t)\big], \; dY_t=\lambda(t) dt,\\
h_t:&=& 
\sum_{j=1}^r \int_0^t\big[ B_s^j(u^{(1)}_s)-B_s^j(u^{(2)}_s)\big]\,dW^j_s.
\end{eqnarray*}
Notice that  almost surely 
\begin{eqnarray*}
\Big| \int_0^T|\delta_t|_V^p\; dY_t\Big|& \leq & 
2^{p-1}\sum_{i=1}^2\int_0^T|u^{(i)}_t|_V^p\, \lambda(t)\,dt<\infty \, ,\\ 
\Big|\int_0^T|z^{\ast}_t|_{V^{\ast}}^q\,dY_t\Big| &\leq &
2^{q-1}\sum_{i=1}^2\int_0^T|A_t(u^{(i)}_t)|^{q}_{V^*}\; \lambda^{1-q}(t)\,dt \\
& \leq & 2^{q-1}\sum_{i=1}^2
\int_0^T\alpha\, |u^{(i)}_t|_V^{p}\; \lambda(t)\,dt 
+  2^{q}\, \int_0^T K_2(t)\,dt<\infty\, ,  
\end{eqnarray*}
and hence almost surely 
\begin{eqnarray*}
\int_0^T|\delta_t|_V\, |z^{\ast}_t|_{V^{\ast}}\,dY_t & \leq & 
\frac{2^{p-1}}{p}\sum_{i=1}^2 \int_0^T|u^{(i)}_t|_V^p\; \lambda(t)\,dt\\
&&  +
\frac{2^{q-1}}{q} \, \alpha \sum_{i=1}^2 \int_0^T |u^{(i)}_t|_V^p\; \lambda(t)\,dt
+  2^{q}\, \int_0^T K_2(t)\,dt<\infty\, . 
\end{eqnarray*}
Thus the conditions of Theorem 1 from \cite{GK:82} on It\^o's
formula holds for the semi-martingale $y$ defined by the
right-hand side of (\ref{ito}).  Hence the monotonicity condition
(C1) yields 
\begin{eqnarray*}
0&\leq&  |\delta_t|^2_H= \int_0^t2\langle \delta_s\, ,\,
z^{\ast}_s\rangle\,dY_s
+[h]_t+m_t\\
&=&2\int_0^t \Big[
\langle u_s^{(1)}-u_s^{(2)},A_s(u_s^{(1)})-A_s(u_s^{(1)})\rangle
+\sum_{j=1}^r |B^j_s(u_s^{(1)})-B^{j}_s(u_s^{(2)})|^2_H\Big]\,ds+m_t
\leq m_t,
\end{eqnarray*}
where $[h]$ is the quadratic variation of $h$, and $m$ is a
continuous local martingale starting from $0$.  By the above
inequality $m$ is non-negative; hence almost surely $m_t=0$ for
all $t\in[0,T]$, which proves that almost surely
$u^{(1)}_t=u^{(2)}_t$ for all $t\in[0,T]$. 

In order to prove the second  statement of the remark we set 
$\gamma(t):=\exp(-\int_0^tK_1(s)\,ds)$, where   
$K_1$ is from condition (C2).  Let $u$ be a solution of 
equation (\ref{u}).  Then by using It\^o's
formula for $\gamma(t)|u(t)|^2_H$ and condition (C2) we get 
$$
\gamma(t)|u(t)|^2_H\leq |u_0|^2_H
+\int_0^t\gamma(s)\bar K_1(s)\,ds+M(t), 
$$ 
where $M$ is a continuous local martingale starting from 0. 
Hence 
$$
E|u(t)|^2_H\leq\gamma^{-1}(T)
\Big[ E|u_0|^2_H+\int_0^T\gamma(s)\bar K_1(s)\,ds\Big] 
$$
for all $t\in[0,T]$, which proves (\ref{becsles}). \hfill $\Box$

\noindent We note that if $u$ is a solution of equation (\ref{u}) 
then under conditions (C2), (C4) and (C5) one can also show 
by standard arguments from \cite{Pa}, \cite{KR} (or see \cite{G})  
that 
$E\big( \sup_{t\in[0,T]}|u_t|^2_H\big) <\infty. $
In the present paper we do not need this estimate, therefore 
we do not prove it.

Our aim is to show that the explicit and implicit
numerical schemes presented below converge to a
stochastic process $u$, which is a solution of
equation (\ref{u}).  Thus, as a byproduct we prove
also the existence part of Theorem \ref{existu}.

First we characterize the solution of equation (\ref{u})
as a minimiser of certain convex functionals.
This characterization, which is a translation of the
method of monotonicity used for example in
\cite{Pa}, \cite{KR} and
\cite{G}, gives a way of proving our approximation theorems.
\smallskip

Fix $T>0$. If  $X$ is a separable Banach space, $\varphi$ is a
positive adapted stochastic process and  $p\in[1,\infty[$, then
${\mathcal L}^p_X(\varphi)$ denotes the Banach space of the
$X$-valued adapted stochastic processes $\{z_t:t\in[0,T]\}$ with
the norm
$$
|z|_{{\mathcal L}^p_X(\varphi)} :=\Big(E \, \int_0^T|z_t|_X^p\,
\varphi(t)\,dt\Big)^{1/p}<\infty\,,
$$
where $|x|_X$ denotes the norm of $x$ in $X$.
If $\varphi=1$, then we use also the notation
${\mathcal L}^p_X$ for ${\mathcal L}^p_X(1)$. 
Let $L_X^p$ denote the Banach space of
$X$-valued random variables $\xi$ with the norm
$$
|\xi|_{L_X^p}:=\big(E|\xi|_X^p\big)^{1/p}.
$$
Let $X$ be embedded in the Banach space $Y$, and let $x=\{x_t\,
:\, t\in[0,T]\}$ and $y=\{y_t\, :\, t\in[0,T]\}$  be stochastic
processes with values in $X$ and $Y$ respectively, such that $x_t(\omega)=y_t(\omega)$ for $dt\times
P$-almost every $(t,\omega)$.  Then we say that $x$  is an $X$-valued
modification of $y$, or that $y$ is a $Y$-valued modification of
$x$.
\begin{Df}                                              \label{A}
Let ${\mathcal A}$ denote the space of triplets $(\xi,a,b)$ satisfying
the following conditions:

\noindent $\bullet$
$\xi : \Omega \rightarrow H$ is
${\mathcal F}_0$-measurable
and such that  $E|\xi|_H^2<+\infty$;

\noindent $\bullet$
$a:[0,T]\times \Omega \rightarrow V^*$ is
a predictable process such that \\
$ E\int_0^T|a_s|^q_{V^*}\:
\lambda^{1-q}(s)\, ds <+\infty$;

\noindent $\bullet$
$b:[0,T]\times \Omega \rightarrow H^r$ is
a predictable process such that\\
$\sum_{j=1}^r
E \int_0^T|b^j_s|^2_H\, ds <+\infty$;

\noindent $\bullet$
There exists a $V$-valued adapted process
$x\in {\mathcal L}^p_V(\lambda)$ such that
\begin{equation}                                                  \label{x}
x_t=\xi + \int_0^t a_s\, ds
+ \sum_{j=1}^r \int_0^t b_s^j\, dW^j_s
\end{equation}
for $dt\times P$-almost all $(t,\omega)\in[0,T]\times\Omega$.
\end{Df}
\smallskip

Let $(\xi,a,b)\in {\mathcal A}$, $x$ defined by (\ref{x}), and
$y\in {\mathcal L}^p_V(\lambda)\cap {\mathcal L}^{2}_H(K_1)$.  Set
\begin{equation}                                                    \label{F}
F_y(\xi,a,b):=E|u_0-\xi|_H^2 +E\int_0^T \Big[ 2\, \langle
x_s-y_s\, ,\, a_s-A_s(y_s)\rangle 
+ \sum_{j=1}^r
|b_s^i-B_s^j(y_s)|^2_H\Big]\, ds\, ,
\end{equation}
and 
 \[ G(\xi,a,b) :=\sup\{ F_y(\xi,a,b)\, :\, y\in {\mathcal
L^p_V(\lambda)}\cap {\mathcal L}^{2}_H(K_1) \}\, .\] 
Due to the growth  condition (C4), for $y\in {\mathcal L}^p(\lambda)$, 
$A_.(y_.)\in {\mathcal L}_{V^*}^{q}(\lambda^{1-q})$. 
Clearly, 
$\langle x\, ,\, z\rangle \in {\mathcal L}^1$ 
for 
$x\in  {\mathcal L}^p(\lambda)$ and 
$z\in {\mathcal L}_{V^*}^{q}(\lambda^{1-q})$, 
by H\"older's inequality.   
Hence (\ref{10.30.11}), (C4) and (C5) imply that 
the functionals $F_y$ and $G$ are well-defined.  
Notice also that $G$ can  take the value
$+\infty$.

\begin{Th}                                                       \label{minimization}
(i) Suppose that conditions (C1)-(C5) hold and 
let $u$ be a solution to (\ref{u}).
Then
\[\inf \{G(\xi,a,b)\, : \, (\xi,a,b)\in {\mathcal A}\}
=G(u_0, A_.(u_.),B_.(u_.))=0\, .\]

(ii)  Assume conditions (C2)-(C5). Suppose that there exist
$(\hat\xi,\hat{a},\hat{b})\in{\mathcal A}$ and some  subset
${\mathcal V}$ of ${\mathcal L}^p_V(\lambda) \cap {\mathcal
L}^{2}_H(K_1)$ dense in ${\mathcal L}^p_V(\lambda)$,  such that 
\begin{equation}                                                    \label{criterion}
F_y(\hat\xi,\hat a,\hat b)\leq 0\; , \quad \forall y\in {\mathcal
V}.
\end{equation}
Then $\hat\xi=u_0$,
\[ u_t=u_0 + \int_0^t \hat{a}_s\, ds
+ \sum_{j=1}^r \int_0^t  \hat{b}_s^j\, dW^j_s, \quad t\in[0,T]\]
is a solution to  (\ref{u}), and $G(u_0,\hat{ a},\hat{b})=0$.
\end{Th}
This theorem,  which is formulated under stronger assumptions in \cite{GM},
is proved in the Appendix for the sake of completeness.
\bigskip

Let $V_n\subset V$ be a finite dimensional subset of $V$
and let $\Pi_n:V^* \rightarrow V_n$
be a bounded linear operator for every integer $n\geq 1$.
Suppose that the following conditions hold:
\medskip

{\bf (H1)} The sequence $(V_n\, ,\, n\geq 1)$ is increasing,
i.e., $V_n \subset V_{n+1}$,  and  $\cup_n V_n$ is dense in $V$.
\medskip

{\bf (H2)} For $x\in V_n$, $\Pi_n x=x$ and for every
$h,k \in H$, $x\in V$ and $y\in
V^*$
\[ (\Pi_n h\, ,\, k)=(h\, ,\, \Pi_n k)\,
\quad \mbox{\rm and}\quad \langle \Pi_n x\, ,\, y\rangle
=\langle x\, ,\,\Pi_n y \rangle \, .\]

{\bf (H3)} For every $h\in H$, $|\Pi_n h |_H \leq  |h|_H$
and $\lim_n |h-\Pi_n h|_H = 0$.
\medskip

For $v\in V_n$, let $|v|_{V_n}=|v|_V$ denote the restriction of the
$V$-norm to $V_n$, and let $|v|_{H_n}=|v|_H$ denote the restriction
 of the $H$-norm to $V_n$.
We denote by $H_n$ the Hilbert space $V_n$ endowed with
the norm $|\, .\,  |_{H_n}$.
 We have
$V_n=H_n\equiv H_n^*=V_n^*$ as topological spaces, where $V_n^*$
is the dual of $V_n$, and $H_n$ is identified with its dual
$H_n^*$ with the help of the inner product in $H_n$.  The
conditions  (H2) and  (H3) clearly imply that
$\Pi_n\circ\Pi_n=\Pi_n$. In particular, if $\{e_i\in V:
i=1.2....\}$ is a complete orthonormal basis in $H$, then the
spaces $V_n:=\mbox{\rm span}(e_i\, ,\, 1\leq i\leq n)$, and the
projections $\Pi_n$ defined by $\Pi_n y:=\sum_{i=1}^n \langle
e_i,y\rangle\, e_i$ for $y\in V^*$ satisfy (H1)-(H3).
\medskip

We now describe several discretization schemes.
Let $m\geq 1$,  and set $\delta_m:=T\, m^{-1}$,
$t_i:=i\delta_m$ for $0\leq i\leq m$.

\subsection{Explicit space-time discretization}
For  $0\leq i\leq m$, $t\in [t_i,t_{i+1}[$ and $1\leq j\leq r$,
define the operators $\tilde{A}_t^m$ and $\tilde{B}_t^{m,j}$  on
$V$ by:
\begin{eqnarray}
\tilde{A}^m_t(x)&:=&\tilde{A}^m_{t_0}(x)
=\tilde{B}^{m,j}_t(x)=\tilde{B}^{m,j}_{t_0}(x)
=0\; \mbox{ \rm for}\; i=0\, , \nonumber\\
\tilde{A}^{m}_t(x)&:=&\tilde{A}^{m}_{t_i}(x)
=\frac{1}{\delta_m}\int_{t_{i-1}}^{t_i}
A_s(x)\, ds\in V^*\;  \mbox{\rm for }\; 1\leq i\leq m\,
,\qquad  \label{tildeAm} \\
\tilde{B}^{m,j}_t(x)&:=&\tilde{B}^{m,j}_{t_i}(x)
=\frac{1}{\delta_m}\int_{t_{i-1}}^{t_i} B^j_s(x)\, ds \in H \;
\mbox{\rm  for}\;  1\leq i\leq m \, .\qquad \label{tildeBm}
\end{eqnarray}
We define an approximation $u^n_m$ of
$u$ by explicit space-time  discretization
of equation (\ref{u}) as follows:
\begin{eqnarray}\label{explicit1}
u^n_m(t)&:=&u^n_m(t_i)\; \mbox{\rm for }\; t\in ]t_i,t_{i+1}[\,
,\quad   0\leq i\leq m-1\, ,\nonumber\\
u^n_m(t_0)&:=& u^n_m(t_1)= \Pi_nu_0\, ,\nonumber \\
u^n_m(t_{i+1})&:= &u^n_m(t_i) + {\delta_m} \,
\Pi_n\tilde{A}^m_{t_i}\big(u^n_m(t_i)\big)
\\
& &+  \sum_{j=1}^r \Pi_n\tilde{B}^{m,j}_{t_i}\big(u^n_m(t_i)\big)\,
\big(W^j_{t_{i+1}}- W^j_{t_i}\big) , \, 1\leq i\leq
m-1.\nonumber
\end{eqnarray}
Notice that the random variables $u^n_m(t_i)$ are
${\mathcal F}_{t_i}$-measurable and \linebreak 
$\Pi_n\tilde{B}^{m,j}_{t_i}\big(u^n_m(t_i)\big)$ is
independent of $\big(W^j_{t_{i+1}}- W^j_{t_i}\big)$.
For every $n\geq 1$ let ${\mathcal B}_n=(e_k\, ,\, k \in I(n)\, )$
 denote a basis of $V_n$, such that
${\mathcal B}_n \subset {\mathcal B}_{n+1}$, and such that
${\mathcal B}=\cup_n {\mathcal B}_n$ is a complete
orthonormal basis of $H$.  For every $n\geq 1$  set
\begin{equation}                                            \label{Cn}
C_{\mathcal B}(n):=\sum_{k\in I(n)} |e_k|_V^2\, .
\end{equation}
The following theorem establishes the convergence of
$u^n_m$ to a solution $u$ of (\ref{u}),
and hence proves the existence of a solution to the equation (\ref{u}).

\begin{Th}                                                  \label{convexplicit}
Suppose conditions (C1)-(C5) with 
$0<\lambda \leq 1$, $p=2$, and conditions (H1)-(H3).
Assume that $n$ and $m$ converge to $\infty$ such that
\begin{equation}                                              \label{link}
\frac{C_{\mathcal B}(n)}{m}\rightarrow 0\,.
\end{equation}
Then the sequence of processes $u^n_m$ converges weakly in
${\mathcal L}^2_V(\lambda)$ to the solution $u$ of equation
(\ref{u}), and $u^n_m(T)$ converges to $u_T$ strongly in $L^2_H$.
\end{Th}
When $D=]0,1[$, $V=W^{1,2}_0(D)$, $H=L^2(D)$, 
$\displaystyle Au=\frac{\partial^2 u}{\partial
x^2}$, and $V_n$ corresponds to the piecewise linear finite elements
methods then condition (\ref{link}) reads $\frac{n^3}{m}\rightarrow
0$.  In this case condition  (\ref{link}) can be weakened substantially. 
(See, e.g., \cite{Gy}).  

\subsection{Implicit discretization
schemes} For every $j=1,\, \cdots ,\, r$ and $i=0,\, \cdots,\,
m-1$ let $A^m$ denote the following average:
\begin{equation}\label{Am}
A^m_t(x):=A^m_{t_i}(x)= \frac{1}{\delta_m} \, \int_{t_i}^{t_{i+1}}
A_s(x)\, ds \; \mbox{\rm for}\; t_i\leq t<t_{i+1}\, .
\end{equation}

We define an approximation $u^m$ for $u$
by an implicit time discretization of equation (\ref{u})
as follows:
\begin{eqnarray}                                                 \label{time2}
u^m(t_0)&:=& 0\, ,\nonumber \\
u^m(t_1)&:=&u_0
+ {\delta_m}\,  A^m_{t_0}\big(u^m(t_{1})\big)\, , \nonumber\\
u^m(t_{i+1})&:= &u^m(t_i)
+ {\delta_m}\, A^m_{t_i}\big(u^m(t_{i+1})\big)\nonumber\\
&&+\sum_{j=1}^r  \tilde{B}^{m,j}_{t_i}\big(u^m(t_i)\big)\,
\big(W^j_{t_{i+1}}- W^j_{t_i}\big)\, ,\quad 1\leq i<m\, ,\nonumber \\
u^m(t)&:=&u^m(t_i)\; \mbox{\rm for }\; t\in ]t_i,t_{i+1}[\, ,\quad
0\leq i<m\, , 
\end{eqnarray}
where the operators $A^m_s$ and $\tilde{B}^{m,j}_s$
have been defined in (\ref{Am}) and (\ref{tildeBm}).
\smallskip

 From the above scheme we get another approximation
$u^{n,m}$ for $u$ by space discretization:

\begin{eqnarray}                                           \label{spacetime2}
u^{n,m}(t_0)&: =&0\, ,                      \nonumber \\
u^{n,m}(t_1)&: =&\Pi_n u_0
+ {\delta_m}\,\Pi_n A^m_{t_0}\big(u^{n,m}(t_{1})\big)\, ,  \nonumber\\
u^{n,m}(t_{i+1})&: = &u^{n,m}(t_i) + {\delta_m} \,
\Pi_n A^m_{t_i}\big(u^{n,m}(t_{i+1})\big)                    \nonumber\\
&&+ \sum_{j=1}^r \Pi_n
\tilde{B}^{m,j}_{t_i}\big(u^{n,m}(t_i)\big)\,
\big(W^j_{t_{i+1}}- W^j_{t_i}\big)\, ,
\quad 1\leq i<m\, ,\nonumber \\
u^{n,m}(t)&: =&u^{n,m}(t_i)\;
\mbox{\rm for }\; t\in ]t_i,t_{i+1}[\, ,\quad
0\leq i<m\, . 
\end{eqnarray}
The following theorem establishes the existence and uniqueness
of $u^m$ and of $u^{n,m}$ for $m$ large enough.

\begin{Th}                                                          \label{existsnm}
Let $p\in [2,+\infty[$ and  assume (C1)-(C5).
Then for  any  sufficiently
large integer $m$ equation (\ref{time2}) has a unique solution
$\{u^m(t_i):i=0,1,...,m\}$ such that $E\big(
|u^{m}(t_i)|_V^p)<+\infty$ for each $i=0,\, \cdots,\, m$.
 If in addition to (C1)-(C5)  conditions (H2) and (H3) also hold,
then there is an integer $m_0\geq 1$ such that for every $m\geq m_0$ and $n\geq1$
equation (\ref{spacetime2}) has a unique solution
$\{u^{n,m}(t_i):i=0,1,...,m\}$ satisfying $E\big(
|u^{n,m}(t_i)|_V^p)<+\infty$ for each $i=0,1,2,...,m$ and $n\geq 1$.
\end{Th}

Once the existence of the solutions to
(\ref{time2}) and to (\ref{spacetime2}) is established,
it is easy to see that $u^m=\{u^m(t):t\in[0,T]\}$ and
$u^{n,m}=\{u^{n,m}(t):t\in[0,T]\}$ are $V$-valued
adapted processes.  Now we formulate our convergence result 
for the above implicit schemes.

\begin{Th}                                               \label{convergence}
 Let $p\in [2,+\infty [$ and assume conditions (C1)-(C5).  Then
for $m\rightarrow \infty$ the sequence of processes $u^m$
converges weakly in  ${\mathcal L}^p_V(\lambda)$ to the solution
$u$ of equation (\ref{u}), and the sequence of random variables
$u^m(T)$ converges strongly  to $u_T$ in $L^2_H$.  If in addition
to (C1)-(C5) conditions (H1)-(H3) also hold, then as $m,n$
converge to infinity, $u^{n,m}$ converge weakly  to the solution
$u$ of equation (\ref{u}) in ${\mathcal L}^p_V(\lambda)$, and the
random variables $u^{n,m}(T)$ converge to $u_T$ strongly in
$L^2_H$.
\end{Th}

\section{Proof of the results}

\subsection{Convergence of the explicit scheme}
We reformulate the equation (\ref{explicit1}) 
in an integral form.
For fixed integer $m\geq 1$ set $t_i:=i\delta_m$,
\begin{equation}
                                                 \label{kappa1}
\kappa_1(t):=t_i
\mbox{ \rm{ for }} \, t\in [t_i,t_{i+1}[, \; 
\mbox{\rm and }\;
\kappa_2(t):=t_{i+1}\,
\mbox{  \rm{  for }} \, t\in ]t_i,t_{i+1}]
\end{equation}
for integers $i\geq0$ and let $\kappa_2(t_0)=t_0$.  Then
(\ref{explicit1}) can be reformulated as follows:
\begin{eqnarray}                                  \label{explicit}
u^n_m(t)&=&\Pi_n u_0 +\int_0^{(\kappa_1(t)-\delta_m))^+}
\Pi_nA_s\big(u^n_m(\kappa_2(s))\big)\, ds\nonumber  \\
&&+ \sum_{j=1}^r
\int_0^{\kappa_1(t)}
\Pi_n\tilde{B}^{m,j}_s\big(u^n_m(\kappa_1(s))\big)\, dW^j_s\, .
\end{eqnarray}

The following lemma provides important bounds for the
approximations.  Set
\[\rho:=\rho(n,m):={\alpha}\, C_{\mathcal B}(n)\delta_m\, ,\]
 and for every $ \gamma\in ]0,1[$,    let
\[ I_{\gamma}=\{ (n,m)\, :\,  n,m\geq1\, ,\;
\rho(n,m)\leq\gamma\}\, ,\]
where $\alpha$ is the constant from
condition (C4), and $C_{\mathcal B}(n)$ is defined by (\ref{Cn}).

\begin{lem}                                                     \label{aprioriexplicit}
Let $p=2$ and conditions (C1)-(C5)  
with $0<\lambda \leq 1$ and (H1)-(H3) hold.  Then for
every $\gamma\in(0,1)$
\begin{eqnarray}
&& \sup_{(n,m)\in I_{\gamma}}\sup_{s\in [0,T]}
E \big| u^n_m(s)\big|_H^2<\infty\, ,                              \label{estimHexplicit}\\
&& \sup_{(n,m)\in I_{\gamma}}\,
E\int_0^T\big| u^n_m(\kappa_2(s))\big|^2_V \,\lambda(s)\,
ds<\infty\,,                                                      \label{estimVexplicit}\\
&& \sup_{(n,m)\in I_{\gamma}}\,
E\int_0^T\big|A_s\big(u^n_m(\kappa_2(s)\big)\big|^2_{V^*}\lambda^{-1}(s) \,
\,  ds <\infty\, ,                                         \label{estimAexplicit}\\
&& \sup_{(n,m)\in I_{\gamma}}\,
\sum_{j=1}^r E\int_0^T \big|  \Pi_n\, \tilde{B}_s^{m,j}
\big(u^n_m(\kappa_1(s)\big)\big|_{H}^2\, ds  <\infty\, .               \label{estimBexplicit}
\end{eqnarray}
\end{lem}

\begin{proof} For any $i=1,\, \cdots,\, m-1$,
\begin{eqnarray*}
E|\,u^n_m(t_{i+1})\, |_H^2&=&E|\,u^n_m(t_i)\,|_H^2
+\delta_m^2\,
E|\Pi_n\tilde{A}_{t_i}^m\big(u^n_m(t_i)\big)|_H^2\\
&& +{\delta_m}\, E\Big[ 2\, \langle u^n_m (t_i), \Pi_n\,
\tilde{A}_{t_i}^m\big( u^n_m (t_i)\big)\rangle 
 + \sum_{j=1}^r
\big|\Pi_n\,
\tilde{B}^{m,j}_{t_i}\big(u^n_m(t_i)\big)\big|_H^2\Big]\, .
\end{eqnarray*}
Adding these equalities, using (H2) and (\ref{tildeAm})  we deduce
\begin{align*}
 E|\, u^n_m&(t_{i+1})\, |_H^2
=E|\,\Pi_n u_0\, |_H^2
+ {\delta_m} \sum_{k=1}^i
E\int_{t_k}^{t_{k+1}} \,
\big| \Pi_n\, \tilde{A}_{t_k}^m\big( u^n_m(t_k)\big)\big|_H^2\, dt\\
&\; +  \sum_{k=1}^{i}
E\int_{t_{k-1}}^{t_k} 2\, \langle \,  u^n_m(t_k)\, ,\,
A_s \big(u^n_m(t_k)\big) \rangle\, ds 
+ \sum_{k=1}^{i} \sum_{j=1}^r
\int_{t_k}^{t_{k+1}}
E\big|\Pi_n\, \tilde{B}^{m,j}_s
\big(u^n_m(t_k)\big)\big|_H^2\, ds \, . \\
\end{align*}
Property (H3), the coercivity condition (C2) and  the growth
condition (C4) with $0<\lambda\leq 1$  and the 
Bunjakovskii-Schwarz inequality yield for every $i=1,\,\cdots,\, m-1$
\begin{align}
E|u^n_m &(t_{i+1})|_H^2  \leq  E|u_0|_H^2
+ {\delta_m}\sum_{k=1}^i \int_{t_{k-1}}^{t_k} 
\sum_{l\in I(n)} E\big[ \langle A_{s}\big(u^n_m(t_k)\big)\, ,
 \, e_l \rangle^2 \big] \, ds \nonumber \\
& + \int_{0}^{t_i} E \Big[ 2\, \big\langle u^n_m(\kappa_2(s)) ,
A_s\big(u^n_m(\kappa_2(s))\big) \big\rangle  
+\sum_{j=1}^r |\Pi_n\,
B^j_s\big(u^n_m(u(\kappa_2(s))\big)|_H^2\Big]\, ds                 \label{equalexplicit} \\
&\leq   E|u_0|_H^2 +
{\delta_m}\, C_{\mathcal B}(n)\, E\int_0^{t_i}
|A_s\big(u^n_m(\kappa_2(s))\big)|_{V^*}^2\, ds
\nonumber\\*
 &\quad  - E\int_0^{t_i} \lambda(s)\,
|u^n_m(\kappa_2(s))|_V^2 \, ds
 + \int_0^{t_i}   \bar{K}_1(s)\, ds  
 + E\int_0^{t_i}   K_1(s)\,|u^n_m(\kappa_2(s))|^2_H\,  ds
\nonumber\\*
 &\leq  E|u_0|_H^2
 - E\int_0^{t_i}\! \!  \lambda(s) \big(1-\alpha \,{\delta_m}\,
C_{\mathcal B}(n)\big)\,
 |u^n_m(\kappa_2(s))|_V^2 \, ds\,
\nonumber \\*
 &\quad  + \int_0^{t_i}  \big[ \bar{K}_1(s)
+ {\delta_m}\, C_{\mathcal B}(n) {K}_2(s)\big]\, ds 
 + E\int_0^{t_i}   K_1(s)\,|u^n_m(\kappa_2(s))|^2_H\,  ds\, .
\nonumber
\end{align}
Hence
\begin{eqnarray}                                 \label{boundsumexplicit}
&&E|u^n_m (t_{i+1})|_H^2+\varepsilon \int_0^{t_i}
E|u^n_m(\kappa_2(s))|_V^2 \, \lambda(s)\,  ds  \leq  E|u_0|_H^2
                \nonumber\\*
&&\;\: + \int_0^{t_i}   K_1(s)\, E|u^n_m(\kappa_2(s))|^2_H\,  ds +
\int_0^{t_i}  \big[ \bar{K}_1(s) +
\alpha^{-1}\gamma{K}_2(s)\big]\, ds\,\qquad 
\end{eqnarray}
for $i=1,\, \cdots ,\, m-1$ and $(n,m)\in I_{\gamma}$,
where $\varepsilon:=1-\gamma  >0 $.  Therefore, the integrability
of $K_1$, $\bar{K}_1$ and ${K}_2$ yields the existence of some positive
constant $C$, which is independent of $n$ and $m$, and the existence of
some positive constants $\alpha_i^m, 1\leq i\leq m$ with $\sup_m\,
\sum_{i=0}^{m-1}\alpha_i^m  <+\infty$, such that
\[E\big[\big|u^n_m \big( k\, \delta_m \big)\big|_H^2\big]\leq C
+ C\sum_{i=0}^{k-1} \alpha_i^m \, E\big[\big|u^n_m\big( i\,
\delta_m \big)\big|_H^2\big]\,
\]
for all $k\in \{1, \cdots, m\}$ and
$(n,m)\in I_{\gamma}$.
Hence by a discrete version of Gronwall's lemma
\begin{equation}                                       \label{boundinHexplicit}
\sup_{(n,m)\in I_{\gamma}}\, \sup_{0\leq i\leq m}
E\big[\big|u^n_m\big( i\, \delta_m \big)\big|_H^2\big] =:
C_{\gamma,\varepsilon} <+\infty\, ,
\end{equation}
which gives (\ref{estimHexplicit}).  The inequalities
(\ref{boundsumexplicit}) with $i=m$ and (\ref{boundinHexplicit})
yield (\ref{estimVexplicit}).  Finally, by (C4), (\ref{10.30.11}),
(\ref{tildeBm}) and (H3) we have:
\[ E\int_0^T\big| A_s\big(u^n_m(\kappa_2(s))\big)\big|_{V^*}^2\,
\lambda^{-1}(s)ds 
\leq \alpha E\int_0^T
|u^n_m(\kappa_2(s))|_V^2\, \lambda(s) \, ds
+\int_0^T{K}_2(s) \, ds,
\]
and for $j=1,\, \cdots,\, r$:
\begin{align*}
E\int_0^T\!|\Pi_{n}\,& \tilde{ B}_t^{m,j}
\big(u^n_m(\kappa_1(t))\big)|_H^2\,dt \leq
\int_0^T \!\! \frac{1}{\delta_m}
\int_{(\kappa_1(t)-\delta_m)^+}^{\kappa_1(t)}\!
E\big|B^j_s\big(u^n_m(\kappa_2(s))\big)\big|_H^2\,ds\,   dt\\
& \leq
E\int_0^T\big|B_s^j\big(u^n_m(\kappa_2(s))\big)|_H^2\, ds 
\leq (2\alpha+1)\,  E\int_0^T\lambda(s)\, |u^{n,m}(\kappa_2(s))|_{V}^p\,ds \\
&   +\int_0^T K_1(s)\, E|u^n_m(\kappa_2(s))|^2_H\,ds
+\int_0^T K_3(s)\,ds  \, .
\end{align*}
Hence (\ref{estimHexplicit}) and (\ref{estimVexplicit}) imply
(\ref{estimAexplicit}) and (\ref{estimBexplicit}).
\end{proof}
\smallskip

\begin{prop}                                       \label{uinftyexplicit}
 Let $p=2$ and conditions (C1)-(C5) 
with $0<\lambda \leq 1$ and (H1)-(H3) hold.  Let
$(n,m)$ be a sequence from $I_{\gamma}$  for some $\gamma\in(0,1)$,
such that $m$ and $n$  converge to infinity.  Then it contains a 
subsequence, denoted also by $(n,m)$, such that:

(i)  $u^{n}_{m}(T)$  converges weakly in $L^2_H$ to some random variable
$u_{T\infty}$,

(ii) $u^{n}_{m}(\kappa_2(\cdot))$  converges weakly in
 ${\mathcal L}^2_V(\lambda)$ to some process  $v_{\infty}$,

(iii) $ A_{\cdot} (u^{n}_{m}(\kappa_2(\cdot)))$
converges weakly in ${\mathcal L}^2_{V^{\ast}}(\lambda^{-1})$
to some process $a_{\infty}$,

(iv)  for any  $j=1,\, \cdots,\, r$,
 $\Pi_{n}\tilde{B}_{\cdot}^{m,j}(u^{n}_{m}(\kappa_1(\cdot))$
converges weakly in ${\mathcal L}^2_{H}$ to some process  $ b^j_\infty$, 

(v) $(u_0,a_{\infty},b_{\infty})\in\mathcal A$,
and for  $dt\times P$-almost every $(t,\omega)\in[0,T]\times\Omega$
\begin{eqnarray}                                      \label{evoluinftyexplicit}
v_\infty(t)&=&u_0+
  \int_0^t a_\infty(s)\,ds
  + \sum_{j=1}^r \int_0^tb^j_\infty(s)\, dW^j(s), \\
                                                 \label{16.29.11}
u_{T\infty}&=&u_0 +
  \int_0^Ta_\infty(s)\, ds
  + \sum_{j=1}^r \int_0^Tb^j_\infty(s)\, dW^j(s) \quad({\mbox{\rm a.s.}}).
\end{eqnarray}
  \end{prop}

\begin{proof}
Assertions (i)-(iv) follow immediately from Lemma
\ref{aprioriexplicit}.  It remains to prove
(\ref{evoluinftyexplicit}) and (\ref{16.29.11}).  Fix $N\geq 1$ and
let $\varphi=\{\varphi(t):t\in[0,T]\}$ be an adapted $V_N$-valued
process such that $|\varphi(t)|_V\leq N$ for all $(t,\omega)$.
From (\ref{explicit})  and (H2), for $n\geq N$ we have
\begin{equation}                                        \label{12.30.11}
 E \int_0^T  \langle u^{n}_{m}(t)\, ,\, \varphi(t)\rangle \,
\lambda(t) \, dt =
 E\int_0^T(u_0, \varphi(t))  \,
\lambda(t)\, dt 
 +J_1+J_2-R-\sum_{j=1}^rR_j,
\end{equation}
with
\begin{eqnarray*}
J_1&:=&E\int_0^T \Big\langle
\int_0^t A_s\big(u^{n}_{m}(\kappa_2(s))\big)\,ds,\varphi(t)
\Big\rangle\,  \lambda(t)\, dt\, ,\\
J_2&:=&\sum_{j=1}^r
E\int_0^T\Big( \int_0^t
\Pi_{n}\, \tilde{B}^{m,j}_s\big(u^{n}_{m}(\kappa_1(s))\big)\, dW^j_s,\,
\varphi(t)\Big) \,
\lambda(t)\, dt\, ,\\
R&:=&E\int_0^T\Big\langle
\int_{(\kappa_1(t)-{\delta_m})^+}^t
A_s\big(u^{n}_{m}(\kappa_2(s))\big)\, ds\,  ,\,
 \varphi(t) \Big\rangle\,  \lambda(t)\,  dt \, ,\\
R_j&:=&E\int_0^T \Big(\int_{\kappa_1(t)}^t
 \Pi_{n}\, \tilde{B}^{m,j}_s\big(u^{n}_{m}(\kappa_1(s))\big)\,
  dW^j_s \, ,\,  \varphi(t)\Big)\,   \lambda(t)\, dt\, .
\end{eqnarray*}
For $(n,m)\in I_{\gamma}$ and $(n,m)\to\infty$, using
(\ref{estimAexplicit}) we obtain
\begin{eqnarray}                                     \label{13.30.11}
&&|R|\leq  N\, E\int_0^T\int_{(\kappa_1(t)-{\delta_m})^+}^t
|A_s\big(u^{n}_{m}(\kappa_2(s))\big)|_{V^{\ast}}\, ds\,dt \nonumber\\
&&\quad \leq  2\, N\, \delta_m\, \Big( E\int_0^T
|A_s(u^{n}_{m}\big(\kappa_2(s))\big)|_{V^{\ast}}^2\, \lambda(s)^{-1} \, ds
\Big)^{\frac{1}{2}}\, T^{\frac{1}{2}} \to 0\, .\qquad 
\end{eqnarray}
 For $j=1,\, \cdots,\, r$ Schwarz's inequality with respect to
$dt\times P$, the isometry of stochastic integrals,
(\ref{estimBexplicit}) and $|\varphi(t)|_H\leq C\,
|\varphi(t)|_V\leq C\, N$ yield:
\begin{eqnarray}                                       \label{14.30.11}
 |R_j| &\leq & C\, 
\left( E\int_0^T |\varphi(t)|_H^2\, dt\right)^{\frac{1}{2}}\;   \left( E\int_0^T 
\left| \int_{\kappa_1(t)}^t \Pi_n\, \tilde{B}^{m,j}_s
\big( u^n_m(\kappa_1(s))\big)\; dW_s^j 
\right|_H^2\, dt \right)^{\frac{1}{2}}\nonumber \\
& \leq & C\, N\, \sqrt{T}\, 
\left( E\int_0^T \int_{\kappa_1(t)}^t \left|
 \Pi_n\, \tilde{B}^{m,j}_s\big( u^n_m(\kappa_1(s))\big)
\right|_H^2\, ds\, dt \right)^{\frac{1}{2}}
\nonumber\\
& \leq & C N\, \sqrt{T\,\delta_m}\,
\left( E\int_0^T\! \! \left|  \Pi_n\, \tilde{B}^{m,j}_s
\big( u^n_m(\kappa_1(s))\big)\right|_H^2\, ds
\right)^{\frac{1}{2}} \rightarrow 0 .
\end{eqnarray}
For $j=1,\, \cdots,\, r$ and $g\in{\mathcal L}^2_H$ let
\begin{equation}
                                                     \label{IS}
F_j(g)(t):=\int_0^t g_s\, dW^j_s,
\quad t\in[0,T]
\end{equation}
 Then by the isometry  of stochastic integrals
\[ \| F_j(g)\|_{{\mathcal L}^2_H(\lambda)}^2= 
 \int_0^T\! \! E\Big( \int_0^t\!  |g_s|_H^2\, ds \Big)\,\lambda(t)\,
dt
 \leq  \int_0^T \lambda(t)\, dt \; \|g\|^2_{{\mathcal L}^2_H}\, ,
\]
which means that the operator $F_j$ defined by (\ref{IS}) is a continuous linear
operator from ${\mathcal L}^2_H$ into ${\mathcal L^2_H(\lambda)}$,
and hence it is continuous also in the weak topologies.   Thus (iv) implies
\begin{equation}                                          \label{15.30.11}
J_2\to \sum_{j=1}^rE\int_0^T\Big(
\int_0^t b^j_{\infty}(s)\,dW_s^j ,\,
\varphi(t)\Big)
\,\lambda(t)\,  dt\, .
\end{equation}
Similarly, the linear operator 
$G: {\mathcal L}_{V^*}^2(\lambda^{-1}) \rightarrow
{\mathcal L}_{V^*}^2(\lambda)$ defined by $G(g)_t=\int_0^t g(s)\, ds$ 
is continuous and hence continuous with respect to the weak topologies. 
Indeed,  
\begin{eqnarray*}
 \|G(g)\|_{{\mathcal L}_{V^*}^2(\lambda)}^2 &\leq &E\int_0^T \lambda(t) 
\Big( \int_0^T \lambda(s)^{-1}\,
|g(s)|^2_{V^*}\, ds \Big)\, \Big(\int_0^t \lambda(s)\, ds\Big) \, dt \\
&\leq& \Big( \int_0^T \lambda(t)\, dt\Big)^2
\, \|g\|^2_{{\mathcal L}^2_{V^*}(\lambda^{-1})}\, .
\end{eqnarray*}
Since $ \varphi \in {\mathcal L}_V^2(\lambda)$,
 (iii) implies
\begin{equation}                                          \label{16.30.11}
J_1\to E\int_0^T \Big\langle \int_0^t a_{\infty}(s)\,ds\, ,\,
\varphi(t)\Big\rangle\,\lambda(t)\, dt .
\end{equation}
Clearly (ii) implies
\begin{equation}                                            \label{17.30.11}
E\int_0^T\big\langle u^{n}_{m}(t),
\varphi(t)\big\rangle \,
\lambda(t) \, dt \to
E\int_0^T (v_\infty(t),\varphi(t))\, \lambda(t) \, dt.
\end{equation}
Thus from (\ref{12.30.11}) we get (\ref{evoluinftyexplicit}) by
(\ref{13.30.11}), (\ref{14.30.11}), (\ref{15.30.11}) -
(\ref{17.30.11}), and by taking into account that $\cup_N V_N$ is
dense in $V$. 
A similar, simpler argument yields that for every random variable
$\psi \in L^2_{V_N}$ such that $E |\psi|_V^2 \leq N$:
\begin{equation}\label{12.30.11bis}
E \left\langle u^n_m(T)\, ,\, \psi \right\rangle = E(u_0\, ,\,
\psi) + \tilde{J}_1 + \tilde{J}_2 - \tilde{R}\, ,
\end{equation}
where as $n,m\rightarrow +\infty$ with $(n,m)\in I_\gamma$,
\begin{eqnarray*}
\tilde{J}_1&=&E \left\langle \int_0^T A_s\big(
u^n_m(\kappa_2(s))\big)\, ds\,
 ,\, \psi \right\rangle\to  
E\left\langle \int_0^T a_\infty(s)\,ds\, ,\, \psi \right\rangle\,  ,\\
\tilde{J}_2&=&\sum_{j=1}^r E \left( \int_0^T \Pi_n
\tilde{B}^{m,j}_s\big( u^n_m(\kappa_1(s))\big)\,dW^j_s\,  ,\, \psi
\right) 
\to  \sum_{j=1}^r E\left( \int_0^T b_\infty^j(s)\,dW^j_s\,
 ,\, \psi \right) \,,\\
|\tilde{R}|&=&E\left| \left( \int_{T-\delta_m}^T A_s\big(
u^n_m(\kappa_1(s))\big)ds\,,\,  \psi \right)\right| \leq C N
\sqrt{\delta_m}\, .
\end{eqnarray*}
Thus, as  $n,m\rightarrow \infty$ with $(n,m)\in I_\gamma$, $E
\left( u^n_m(T)\, ,\, \psi \right) \rightarrow E(u_{T \infty}\,,\,
\psi)$. Since $\cup_N V_N$ is dense  in $H$, this concludes the
proof.
\end{proof}
\smallskip

\begin{prop}\label{exp}
Let $p=2$, (C1)-(C5) with $0<\lambda \leq 1$ 
 and (H1)-(H3) hold. Let $(n,m)$ be a
sequence of pair of positive integers such that 
$m$ and $n$ converge to infinity, and  $C_{\mathcal
B}(n)/{m}\rightarrow 0$.  Then the assertions of Proposition
\ref{uinftyexplicit} hold, and for every $y\in{\mathcal L}^p(\lambda)$:
\begin{equation}                                   \label{ineqyexp}
\int_0^T E\Big[ 2\langle v_\infty(t) - y(t)\, , \,
a_\infty (t)-A_t(y(t))\rangle + \sum_{j=1}^r
|b^j_\infty (t)-B^j_t(y_t)|_H^2 \Big]\, dt \leq 0\, .
\end{equation}
The process $v_{\infty}$ has an $H$-valued continuous
modification, $u_{\infty}$, which is the solution of equation
(\ref{u}), and $E|u^n_m(T)-u_\infty(T)|_H^2\to 0$.
\end{prop}

\begin{proof} Since $C_{\mathcal B}(n)/{m}\rightarrow 0$,
with finitely many exceptions all pairs $(n,m)$ from the given
sequence belong to $I_{\gamma}$.  Thus we can apply Proposition
\ref{uinftyexplicit} and get a subsequence, denoted again by
$(n,m)$, such that assertions (i)--(v) of Proposition
\ref{uinftyexplicit} hold. Notice that $v_{\infty}\in{\mathcal
L}_V^2(\lambda)$ and $a_{\infty}\in{\mathcal
L}_{V^{\ast}}^2(\lambda^{-1})$.  Thus from
(\ref{evoluinftyexplicit}) by Theorem 1 from \cite{GK:82} on
It\^o's formula we get that $v_{\infty}$ has an $H$-valued
continuous modification $u_{\infty}$, and a.s.
\begin{equation}                                    \label{22.30.11}
E|u_{\infty}(T)|_H^2
=E|u_0|^2_H+E\int_0^T\Big[ 2\, \langle v_{\infty}(s),a_{\infty}(s)\rangle
+\sum_{j=1}^r|b^j_{\infty}(s)|_H^2\Big]\,ds\, . 
\end{equation}
Moreover, by (\ref{evoluinftyexplicit}) and
(\ref{16.29.11}) we get $u_{\infty}(T)=u_{T\infty}$.  For
$y\in{\mathcal L}^2_V(\lambda)$ such that \linebreak
 $\sup_{0 \leq
t\leq T} E|y_t|_H^2 <+\infty$, let
\begin{eqnarray*}
 F^n_m(y)
&:=& E\int_0^T\Big[ 2\, \big\langle u^{n}_{m}(\kappa_2(t))-y(t)\,
,\,  A_t(u^{n}_{m}\big(\kappa_2(t))\big) -
A_t(y_t)\big\rangle \\ &&
 +\sum_{j=1}^r  \big|
\Pi_{n}B^j_t(u^{n}_{m}\big(\kappa_2(t))\big) -
\Pi_{n}B^j_t(y(t))\big|_H^2\Big]\,dt\, .
\end{eqnarray*}
 Notice that the growth condition (C4) and (\ref{10.30.11}) imply
that for $x,z\in {\mathcal L}^2_V(\lambda)$, 
$\langle x_.\, ,\, A_.(z_.)\rangle \in {\mathcal L}^1$ and 
 $B^j(y)\in {\mathcal L}^2_H(K_1)$ for $j=1,\, \cdots,\, r$;  since the
estimates (\ref{estimVexplicit}), (\ref{estimAexplicit}) and
(\ref{10.30.11}) hold, $F^n_m(y)$ is well-defined and  is finite.
By the monotonicity condition (C1), (H3) and by inequality
(\ref{equalexplicit})
\begin{equation}
0\geq F^n_m(y) \geq E|u^{n}_{m}(T)|_H^2 - E|u_0|_H^2 +2\,
  E\int_0^T\langle y_t \, ,\, A_t(y_t)\rangle \, dt
-R-2J_1-2J_2-2J_3+J_4,  \label{ineqykex}
\end{equation}
with
\begin{eqnarray*}
&&R:=\delta_mE\int_0^{T-\delta_m}
\sum_{l\in I(n)}\langle A_s\big(u^{n}_{m}(\kappa_2(s))\big)\, ,\, e_l \rangle^2\, ds\, , \\
&&J_1:=E\int_0^T\langle u^{n}_{m}(\kappa_2(t))\,,
\, A_t(y_t)\rangle\, dt\, ,\\
&&J_2:=E\int_0^T\left\langle y_t\,,\,
A_t\big(u^{n}_{m}(\kappa_2(t))\big)\right\rangle\, dt\, ,  \\
&&J_3:=\sum_{j=1}^r  E\int_0^T
    \left( \Pi_{n}B_t^j\big(u^{n}_{m}(\kappa_2(t))\big)\, ,\,
   B_t^j(y_t) \right)\, dt\, ,\\
&& J_4:=\sum_{j=1}^r E \int_0^T |\Pi_{n}B_t^j(y_t)|_H^2\, dt.
\end{eqnarray*}
Since $\lambda^{-1}\geq 1$, (\ref{estimAexplicit}) 
implies that for $C_{\mathcal B}(n)/m\to 0$
\begin{equation}                                                      \label{18.30.11}
|R|\leq T\, \frac{C_{\mathcal B}(n)}{m}\, E
\int_0^T \big|A_s(u^n_m(\kappa_2(s)))\big|_{V^{\ast}}^2\,ds  \to 0.
\end{equation}
By Proposition \ref{uinftyexplicit}, as $C_{\mathcal B}(n)/m\to 0$,
\begin{eqnarray}
J_1&=&E\int_0^T\big\langle u^{n}_{m}(\kappa_2(t))\, ,\, A_t(y_t)\,
\lambda(t)^{-1}\big\rangle \,\lambda(t)\,  dt 
\rightarrow
E\int_0^T\big\langle u_\infty(t) \, ,\, A_t(y_t)\, \big\rangle\, dt\, ,   \label{19.30.11}\\
J_2&=&E\int_0^T\big\langle \lambda(t)y_t\, , \,
A_t\big(u^{n}_{m}(\kappa_2(t))\big) \big\rangle\,
\lambda^{-1}(t)dt  
\rightarrow  E\int_0^T \big\langle y_t\, , \,
a_\infty(t)\, \big\rangle \, dt. \label{20.30.11}
\end{eqnarray}
Notice that
\[
 E\int_0^T\Big(\Pi_{n}\tilde B_t^{m,j}\big(u^{n}_{m}(\kappa_1(t))\big)
\, ,\,  B_t^j(y_t)) \Big)\,dt
  =E\int_0^T\Big(\Pi_{n}B_t^j\big(u^{n}_{m}(\kappa_2(t))\big)
\, ,\,  S_mB_t^j(y_t)) \Big)\,dt\, ,
\]
where $S_m$ is the averaging operator, defined by
\begin{equation}
(S_mZ)(t):=\left\{
\begin{array}{lll}
\delta_m^{-1}\int_{\kappa_1(t)+{\delta_m}}^{\kappa_1(t)+2{\delta_m}} Z_s\,
ds &\mbox{\rm if}& 0\leq t\leq T-\delta_m\, ,\\
0&\mbox{\rm if}& T-\delta_m<t\leq T \,
\end{array}\right.                                            \label{moyennebis}
\end{equation}
for $Z\in{\mathcal L}^2_H$.
Hence, taking into account Proposition \ref{uinftyexplicit} (iv) and
$$
\lim_{m\to\infty}E \int_0^T|(S_mZ)_t-Z_t|_H^2\,dt=0\; , \quad
\forall Z\in{\mathcal L}^2_H\, ,
$$
as $C_{\mathcal B}(n)/m\to 0$  we get
\begin{equation}                                                     \label{17.30.11bis}
J_3\to \sum_{j=1}^r E \int_0^T \Big(  b^j_\infty(t),
B_t^j(y_t)\Big)\, dt\, .
\end{equation}
Using (H3) and the dominated convergence theorem, since
$B_.(y_.)\in {\mathcal L}^2_H$,
we obtain 
\begin{equation}      \label{17.30.11ter}
J_4 \to \sum_{j=1}^r E\int_0^T |B_t(y_t)|_H^2\, dt\, .
\end{equation}
Since $u^{n}_{m}(T)$ converges weakly in $L^2_H$ to $u_{T\infty}=u_{\infty}(T)$,
\begin{equation}                                                   \label{21.30.11}
d := \liminf_{n,m\to\infty} E|u^{n}_{m}(T)|_H^2
-E|u_\infty(T)|_H^2 \geq 0\,.
\end{equation}
Letting $n,m\to\infty$ with $C_{\mathcal B}(n)/m\to 0$ in
(\ref{ineqykex}) and using  (\ref{22.30.11}), (\ref{18.30.11}) -
 (\ref{20.30.11}) and
 (\ref{17.30.11bis}) - (\ref{21.30.11}), we deduce that for $y\in {\mathcal L}^2_V(\lambda)$
 with $\sup_t E|y_t|_H^2<+\infty$ and $F_y$ defined by (\ref{F}):
\begin{eqnarray}
0&\geq&\;  d + E|u_\infty(T)|_H^2  - E|u_0|_H^2  + 2\,
E\int_0^T\langle y_t\, ,\, A_t(y_t)\rangle\, dt 
 -2\, E\int_0^T\langle u_\infty(t)\, ,\,
A_t(y_t)\rangle\, dt  \nonumber \\
& &\quad
 -2 E\int_0^T\langle y_t\, ,\,  a_\infty(t) \rangle\, dt 
+ \sum_{j=1}^r  E\int_0^T\Big[ |B_t^j(y_t)|_H^2 
-2\, (b_\infty^j(t), B_t^j(y_t))\Big]\, dt
\nonumber \\
&=&\; d+F_y(u_0,a_{\infty},b_{\infty}),  \label{inegalex}
\end{eqnarray}
by (\ref{18.30.11}) - (\ref{20.30.11}), 
(\ref{17.30.11bis}) - (\ref{21.30.11}), and
taking into account (\ref{22.30.11}). 
Hence by Theorem \ref{minimization} (ii), $u$ is a solution to
equation (\ref{u}).  Taking $y:=u$ in the above inequality we get
$d\leq 0$, and hence $d=0$.  Thus the approximations $u^{n}_{m}(T)$
converge weakly in $ L^2_H$ and their $ L^2_H$-norms converge to
that of $u_\infty(T)$, which imply the strong convergence of
$u^{n}_{m}(T)$ in  $ L^2_H$ to $u(T)$.
\end{proof}
\smallskip

Now we  conclude the proof of Theorem \ref{convexplicit}. 
Let  $(n,m)$ be a sequence of pairs of positive integers such
that  $m$ and $n$ converge to infinity and 
 $C_{\mathcal B}(n)/{m}\rightarrow 0$;  the previous
proposition proves the existence of  a subsequence of the explicit
approximations $u^n_m$ that converges weakly in 
${\mathcal L}^2_V(\lambda)$ to a solution $u_{\infty}$ of the equation
(\ref{u}), and such that  $u^n_m(T)$ converges strongly to
$u_{\infty}(T)$ along the same subsequence.  Since by  Remark
\ref{uniqueness} the solution to (\ref{u}) is unique, the whole
sequence $u^n_m$ converges weakly in ${\mathcal L}^2_V(\lambda)$
to the solution of the equation (\ref{u}), and the whole sequence
$u^n_m(T)$ converges strongly in
 $L^2_H$  to  $u_{\infty}(T)$. \hfill{$\Box$}

\subsection{Existence and uniqueness of solutions to the implicit schemes}
The following proposition establishes existence and uniqueness
of the solution to the equation  $Dx=y$  and provides
an estimate of the norm of $x$ in terms of that of $y$.
\begin{prop}\label{infdim}
Let  $D:V\rightarrow V^*$ be such that:

(i) $D$ is monotone, i.e., for every $x,y\in V$,
$\langle D(x) -D(y),x-y\rangle  \geq 0$.

(ii) $D$ is hemicontinuous, i.e., 
${\displaystyle \lim_{\varepsilon\rightarrow 0}
  \langle D(x+\varepsilon y),z\rangle  =\langle D(x),z\rangle }$
for every $x,y,z\in V$. 

(iii) $D$ satisfies the growth condition, i.e.,
there exists $K>0$ such that for every $x\in V$,
\begin{equation}                                    \label{11.04.11}
|D(x)|_{V^*}\leq K\, (1+|x|_V^{p-1}).
\end{equation}

(iv) $D$ is coercive, i.e., there exist constants
$C_1>0$ and $C_2\geq 0$ such that
\[ \langle D(x),x\rangle \geq C_1\, |x|_V^p - C_2\; , \quad \forall x\in V\, .\]

Then for every $y\in V^*$, there exists $x\in V$ such that $D(x)=y$ and
\begin{equation}                                        \label{norm}
|x|_V^p \leq \frac{C_1 +2\, C_2}{C_1} + \frac{1}{C_1^2}\, |y|_{V^*}^2\, .
\end{equation}

If there exists a positive constant $C_3$ such that
\begin{equation}
\label{strongmono}
 \langle D(x_1)-D(x_2),x_1-x_2\rangle \geq C_3\, |x_1-x_2|^2_{V^*}\; ,
\quad \forall x_1,x_2\in V\, ,
\end{equation}
then for any $y\in V^*$, the equation $D(x)=y$ has a unique solution $x\in V$.
\end{prop}
This result is known, or can easily be obtained from
well-known results.  We include its proof in the Appendix for the convenience of the reader.
\bigskip

{\it Proof of Theorem \ref{existsnm}.}
To  prove this theorem, we need to  check the conditions of the previous proposition
for the operators $D: V\to V^{\ast}$ and $D_n:V_n\to V_n$, defined by
\[ D:=I-\int_{t_i}^{t_{i+1}} A_s\, ds \quad \mbox{\rm and }\quad
D_n:=I_n-\int_{t_i}^{t_{i+1}} \Pi_nA_s\, ds \] for each
$i=0,1,2,...,m-1$, where  $I : V\rightarrow V^*$ denotes the canonical embedding and
  $I_n$ denotes the identity operators  on $V_n$.
Hence $u^{m}(t_i)$ and $u^{n,m}(t_i)$ can be uniquely defined recursively for
 $0\leq i\leq m$ by the equations (\ref{time2}) and (\ref{spacetime2}), respectively.

We first check that $D$ satisfies the strong monotonicity
condition.  Let $x,y\in V$.  Then (C1) implies
\[ \langle D(x)-D(y),x-y\rangle = |x-y|_H^2
- \int_{t_i}^{t_{i+1}} \langle A_s(x)-A_s(y), x-y\rangle \, ds
\geq |x-y|_H^2\, . \]
Let us check that $D$ is hemicontinuous.  Let
$x,y,z\in V$  and $\varepsilon \in \mathbb{R}$:
\[ \langle D(x+\varepsilon y),z\rangle
= \langle x+\varepsilon\,  y,z\rangle - \int_{t_i}^{t_{i+1}}
\langle A_s(x+\varepsilon y) \, ,\, z\rangle\, ds \, .\] As
$\varepsilon \rightarrow 0$, condition (C3) implies that for every
$s\in [t_i, t_{i+1}]$, $\langle A_s(x+\varepsilon y) \, ,\,
z\rangle$ converges to $\langle A_s(x) ,z\rangle$.   Hence, using
condition (C4)  we have that for every $\varepsilon \in ]0,1]$:
\[ |\langle A_s(x+\varepsilon y\, ,\, z\rangle | \leq  C\, \alpha^{\frac{1}{q}}\, \lambda(s)\, 
\Big[ |x|_V^{p-1} + |y|_V^{p-1}\Big]\, |z|_{V^*}
  + \frac{1}{p} \, \lambda(s) + \frac{1}{q}\, K_2(s)
\in {\mathcal L}^1\, .
\]

Thus, we get the hemicontinuity of $D$ by  the Lebesgue
theorem on dominated convergence.

We check that the operator $D$  satisfies the growth condition
(\ref{11.04.11}).  Let $x\in V$.   Then by condition (C4) for $p\in [2,+\infty [$  
we have
\[ |D(x)|_{V^*}\leq  |x|_{V^*} 
+  C_1\int_{t_i}^{t_{i+1}} \Big[ \lambda(s)\,
|x|_V^{p-1} +  \frac{1}{p}\, \lambda(s) 
+ \frac{1}{q}\, K_2(s) \Big] \, ds
\leq C_2\, [1+|x|_V^{p-1}] \,
\]
with some constants $C_1,C_2$. 

We check that for $m$ large enough, $D$  satisfies the coercivity
condition.  Let $x\in V$; then using (C2) we have:
\begin{align*}
\langle D(x),x\rangle& \geq  \; |x|_H^2 + \int_{t_i}^{t_{i+1}} \frac{1}{2}\,
\big[\lambda(s)\, |x|_V^p - K_1(s)\, |x|_H^2
-\bar{K}_1(s) \big]\, ds\\
&\geq\;   \frac{1}{2}\,
\left( \int_{t_i}^{t_{i+1}} \lambda(s)\, ds \right)\, |x|_V^p - \frac{1}{2}\,
\int_{t_i}^{t_{i+1}} \bar{K}_1(s)\, ds
+ |x|_H^2\,
\left[ 1-\frac{1}{2}\, \int_{t_i}^{t_{i+1}} K_1(s)\, ds \right]\, .
\end{align*}
Since $K_1$ is integrable, for large enough $m$,
${\delta_m}=t_{i+1}-t_i = T\, m^{-1}$ is small enough to imply
that $\int_{t_i}^{t_{i+1}} K_1(s)\, ds <2$; thus $D$ is coercive.
Using (H2), (H3) and the equivalence of the norms $|\, .\,
|_{V_n}$, $|\, .\, |_{H_n}$ and $|\, .\, |_{V_n^*}$ on $V_n$,   
similar arguments show that $D_n$ satisfies conditions (i)-(iv) too.
\smallskip

We finally prove by induction that the random variables
$u^{n,m}(t_i)$ and $u^m(t_i)$ belong to $L^p_V$.  This is obvious
for $t_0=0$, and for $t_1$ it follows immediately from the
estimate (\ref{norm}).   Let $i$ be an integer in $\{1,\, \cdots,\,
m-1\}$,  assume that $E(|u^{n,m}(t_k)|_V^p)< +\infty$ for every
$k\in \{1,\, \cdots ,\, i\}$ and set
\[
y= u^{n,m}(t_i)+  \sum_{j=1}^r \int_{t_i}^{t_{i+1}}
\Pi_n\tilde{B}^{m,j}_s(u^{n,m}(t_i))\, dW_s^j \in V_n\,.
\]
Then by the isometry of stochastic integrals, and by Remark
\ref{15.29.10} for $p\in [2, +\infty [$ we have
\begin{eqnarray*}
E|y|_{V^*}^2
&\leq&  C_1\, E|u^{n,m}(t_i)|_V^2+ C_1\, E\int_{t_i}^{t_{i+1}}
|\Pi_n\tilde{B}_s^{m,j} (u^{n,m}(t_i))|_H^2\, dt \\
&\leq& C\Big[1+\int_{t_{i-1} }^{t_i} \lambda(s)\, ds\Big]\,
E|u^{n,m}(t_i)|_V^p\\
&&\; +C\, \int_{t_{i-1}}^{t_i} K_1(s)\,
E|u^{n,m}(t_i)|_H^2\, ds 
+  C\Big[1+\int_{t_{i-1}}^{t_i}K_3(s)\,ds\Big]<\infty\,.
\end{eqnarray*}
Hence (\ref{norm}) shows that $E(|u^{n,m}(t_{i+1})|_V^p)<+\infty$.
In the same way we get the finiteness of the $p$-th moments of the
$V$-norm of $u^m(t_i)$ for $i=0,1,2,..m$.  \hfill $\Box$

\subsection{Convergence of the implicit schemes}

We first prove some a priori estimates on the processes
$u^m$ and $u^{n,m}$ and give an evolution
formulation to the equations satisfied by these processes.

Recall that for $0\leq i<m$ and $t\in ]t_i,t_{i+1}[$, we  set
$\kappa_1(t)=t_i$ and $\kappa_2(t)=t_{i+1}$, while  for  $i=0,
\cdots, m$, we set $\kappa_1(t_i)=\kappa_2(t_i)=t_i$.   Let $A^m$
and $\tilde{B}^{m,j}$ be defined in (\ref{Am}) and
(\ref{tildeBm}).  Then equations (\ref{time2}) and
(\ref{spacetime2}) can be cast in the integral form:
\begin{equation}                                        \label{time}
 u^m(t)=u_0\,  1_{\{t\geq t_1\}}
+\int_0^{\kappa_1(t)}  A_s\big(u^m(\kappa_2(s))\big)\, ds 
+ \sum_{j=1}^r
\int_0^{\kappa_1(t)} \tilde{B}^{m,j}_s\big(u^m(\kappa_1(s))\big)\, dW^j_s\,,
\end{equation}
and
\begin{equation}                                    \label{spacetime}
 u^{n,m}(t)=\Pi_nu_0  1_{\{t\geq t_1\}}  +
\int_0^{\kappa_1(t)}\!\! \Pi_n A_s\big(u^{n,m}(\kappa_2(s))\big)\, ds 
 + \sum_{j=1}^r
\int_0^{\kappa_1(t)}\!\! \Pi_n\tilde{B}^{m,j}_s(u^{n,m}\big(\kappa_1(s))\big)\, dW^j_s,
\end{equation}
respectively.
\begin{lem}                                      \label{apriori}
Let conditions (C1)-(C5) and (H1)-(H3) hold.  Then there exist an
integer $m_1$ and some constants $L_i, 1\leq i\leq 4$, such that:
\begin{eqnarray}
&&\sup_{s\in [0,T]} E\big| u^{n,m}(s)\big|_H^2\leq L_1\,, \label{estimH}\\
&& E\int_0^T \big| u^{n,m}(\kappa_2(s))\big|^p_V \,
\lambda(s)\, ds  \leq L_2\, ,                             \label{estimV}\\
&&E\int_0^T \big| A_s\big(u^{n,m}(\kappa_2(s))\big)
\big|^q_{V^*}\,\lambda(s)^{1-q} \, ds
\leq L_3\, ,                                             \label{estimA}\\
&&\sum_{j=1}^rE\int_0^T 
\big| \tilde{B}^{m,j}_s\big(u^{n,m}(\kappa_1(s))\big)
\big|_{H}^2 \, ds \leq L_4\,                                   \label{estimB}
\end{eqnarray}
for all $m\geq m_1$ and $n\geq1$. 
Under conditions (C1)-(C5) the above estimates hold with the
implicit approximations $u^m$ in place of $u^{n,m}$ for
all sufficiently large $m$.
\end{lem}
\begin{proof}
We only prove the estimates for $u^{n,m}$.  The proof of the
estimates for $u^m$ is essentially the same, and we omit it.
We set $\Delta W^j_{t_i}=W^j_{t_{i+1}}-W^j_{t_i}$
for $i=0, \cdots, m-1$, $j=1\, \cdots, r$.  Then from the definition
of the approximations $u^{n,m}$ we get
\[
|u^{n,m}(t_1)|_H^2-|\Pi_nu_0|^2_H
=2\langle u^{n,m}(t_1), A^m_0(u^{n,m}(t_1))\rangle\delta_m 
-\, |\Pi_n A^m_0(u(t_1))|_H^2\delta_m^2
\]
and for $i=1,\, \cdots, \, m -1$:
\begin{align*}
u^{n,m}(t_{i+1})|_H^2 & - |u^{n,m}(t_{i})|^2_H =2\big\langle
u^{n,m}(t_{i+1}), A_{t_i}^m(u^{n,m}(t_{i+1}))\big\rangle\,
\delta_m 
-\big|\Pi_n
A_{t_i}^m(u^{n,m}(t_{i+1}))\big|_H^2\delta_m^2 \\
& +2\sum_{j=1}^r \big(u^{n,m}(t_{i}), \Pi_n \tilde
B^{m,j}_{t_i}(u^{n,m}(t_{i}))\big)\, \Delta W^j_{t_i} 
 +\Big|\sum_{j=1}^r \Pi_n \tilde
B^{m,j}_{t_i}(u^{n,m}(t_{i})\,\Delta W^j_{t_i}\Big|_H^2 \, .
\end{align*}
Hence adding these equations and taking expectation we obtain
\begin{align*}
 E|u^{n,m}&(t_k)|_H^2= E|\Pi_n(u_0)|_H^2  + 2\,
E\int_0^{t_{k}} \big\langle \,
 A_s \big(u^{n,m}(\kappa_2(s))\big),  
u^{n,m}(\kappa_2(s))\big\rangle\, ds \\
 &  +  \sum_{j=1}^r  E\int_{t_1}^{t_{k}} \big|\Pi_n\tilde{B}^{m,j}_s
\big(u^{n,m}(\kappa_1(s))\big)\big|_H^2\, ds
  - {\delta_m} E\int_0^{t_k}
\big| \Pi_nA_t \big(u^{n,m}(\kappa_2(s))\big)\big|_H^2\, ds\,
\end{align*}
for $k=1,2,...,m$, which implies
\begin{align}
E|u^{n,m}&(t_k)|_H^2\leq E|u_0|_H^2 - {\delta_m}\,
E\int_0^{t_k}\! \big|\Pi_n A_s
\big( u^{n,m}(\kappa_2(s))\big)|_H^2\, ds \nonumber \\*
&\, +E\int_0^{t_k}\!\Big\{
2\, \big\langle A_s\big(u^{n,m}(\kappa_2(s))\big)\, ,\,
u^{n,m}(\kappa_2(s))\big\rangle 
+\sum_{j=1}^r 
\big|B_s^j(u^{n,m}(\kappa_2(s)))\big|_H^2 \Big\}\, ds   \label{equal} \\*   
&\leq \; E|u_0|_H^2 - E\int_0^{t_k} \lambda(s)\,
|u^{n,m}(\kappa_2(s))|_V^p \, ds
 + E\int_0^{t_k} K_1(s) \,
|u^{n,m}(\kappa_2(s))|_H^2\, ds \nonumber  \\*
& \, + \int_0^{t_k} \bar{K}_1(s)\, ds\, ,\nonumber
\end{align}
by the definition of $\tilde{B}^{m,j}$, the coercivity condition
(C2), and by (H2).   For $m$ large enough,
 $\gamma_m=\sup\{ \int_{t_{k-1}}^{t_k} K_1(s)\, ds\; :\; 1\leq k\leq m \} \leq \frac{1}{2}$.
Consequently, there exists an integer $m_1$ such that for all $n\geq 1$, $m\geq m_1$
and $k=1,2,...,m$:
\begin{equation}                                   \label{boundsum}
\frac{1}{2}\, E|u^{n,m}(t_k)|_H^2 + E\int_0^{t_k} \!
|u^{n,m}(\kappa_2(s))|_V^p\, \lambda(s)\,  ds 
\leq C +
\int_0^{t_{k-1}}\!\!   K_1(s)\, E|u^{n,m}(\kappa_2(s))|_H^2\, ds
\, .
\end{equation}
Hence a discrete version of Gronwall's lemma implies the existence of a constant
 $C>0$ such that
\begin{equation}                                    \label{boundinH}
\sup_{n\geq 1}\, \sup_{m\geq m_1}\,
\sup_{0\leq k\leq m}
E\big|u^{n,m}\big(k\, \delta_m\big)\big|_H^2=C<\infty\,,
\end{equation}
which implies (\ref{estimH}).  The inequalities (\ref{boundsum})
and (\ref{boundinH}) yield (\ref{estimV}).
Notice that by the growth condition (C4)
\[ E\int_0^T\big| A_s(u^{n,m}\big(\kappa_2(s))\big)\big|_{V^*}^q\,
\lambda(s)^{1-q}\,  ds 
\leq \alpha\,  E\int_0^T
 |u^{n,m}(\kappa_2(s))|_V^p\,\lambda(s)\,ds + \int_0^T{K}_2(s)\, ds
\]
 and by the definition of $\tilde{B}^{m,j}$ and by Remark
\ref{15.29.10},
\begin{align*}
 E&\int_0^T|\tilde{B}^{m,j}_s(u^{n,m}(\kappa_1(s)))|_H^2\, ds \leq
E \int_0^{T-{\delta_m}} \! \big|  B_s^j (u^{n,m}(\kappa_2(s))\big|_{H}^2 \, ds
\\
& \leq (2 \alpha+1) E\int^T_0 \!\! |u^{n,m}(\kappa_2(s))|_V^p\,\lambda(s)\,ds
+E\int_0^T \!\! K_1(s)|u^{n,m}(\kappa_2(s))|_H^2\,ds+
\int_0^T \!\!  K_3(s)ds.
\end{align*}
Thus estimates (\ref{estimH}) and (\ref{estimV}) imply estimates
(\ref{estimA}) and (\ref{estimB}).
\end{proof}

\begin{prop}\label{uinfty}
Let conditions (C1)--(C5) and (H1)-(H3) hold.
Then for any sequence $(n,m)\to\infty$ of pairs of positive integers
there exists a subsequence, denoted also by (n,m), such that:

(i)  $u^{n,m}(T)$  converges weakly to $u_{\infty T}$ in
$L^2_H$,

(ii) $u^{n,m}(\kappa_2(.))$
converges weakly in ${\mathcal L}^p_V(\lambda)$ to $v_{\infty}$,

(iii)
$A_. (u^{n,m}(\kappa_2(.)))$  converges weakly
in ${\mathcal L}^q_{V^{\ast}}(\lambda^{-1})$ to  $a_{\infty}$,

(iv) $\Pi_{n}\tilde{B}^{m,j}_.(u^{n,m}(\kappa_1(.)))$   converges
weakly in  ${\mathcal L}^2_H$   to $b^j_{\infty}$ for each
$j=1,2,...,r$.

(v) $(u_0,a_{\infty},b_{\infty})\in\mathcal A$,
and for  $dt\times P$-almost every
$(t,\omega)\in[0,T]\times\Omega$
\begin{eqnarray}                              \label{evolimplicit}
v_\infty(t)&=& u_0+
  \int_0^t a_\infty(s)\,ds
  + \sum_{j=1}^r \int_0^tb^j_\infty(s)\, dW^j(s), \\
                                              \label{16.31.11}
u_{T\infty}&=&u_0 +
  \int_0^Ta_\infty(s)\, ds
  + \sum_{j=1}^r \int_0^Tb^j_\infty(s)\, dW^j(s) \quad({\mbox{\rm a.s.}}),
\end{eqnarray}
\begin{equation}                                                  \label{20.01.11}
F_y(u_0,a_{\infty}, b_{\infty})\leq 0\; ,  \quad
\forall y\in{\mathcal L}^p_V(\lambda)\,.
\end{equation}

(vi) The process $v_{\infty}$ has an $H$-valued continuous
modification, $u_{\infty}$, which is the solution of equation
(\ref{u}).   Moreover, the sequence $u^{n,m}(T)$ converges strongly
in $L^2_H$ to  $u_{\infty T}=v_{\infty}(T)$.

Under conditions (C1)-(C5) the above assertions hold with $u^{m}$
and $\tilde{B}^{m,j}$, in place
of  $u^{n,m}$ and $\Pi_n\tilde{B}^{m,j}$, respectively.
\end{prop}

\begin{proof} We prove the lemma for subsequences of $u^{n,m}$.
The proof for the sequence $u^m$ is essentially the same, and we
omit it.   The assertions (i)-(iv) are immediate consequences of
Lemma \ref{apriori}. We need only prove  assertions (v) and (vi). 
For fixed $N\geq 1$ let $\varphi=\{\varphi(t):t\in[0,T]\}$ be a
 $V_N$-valued adapted stochastic process such that $|\varphi(t)|_V\leq N$
for all $t\in[0,T]$ and $\omega\in\Omega$.
Then from equation (\ref{spacetime}) for $n\geq N$ we have
\begin{eqnarray}
&&E\int_0^T( u^{n,m}(t),\varphi(t))\, \lambda(t) \, dt
 =E\int_0^T  1_{\{t\geq t_1\}}\, \big(\Pi_{n}u_0\, ,\,
\varphi(t)\big)  \,  \lambda(t)  \, dt                                \nonumber\\
&&\qquad\quad +  E\int_0^T \Big\langle \int_0^{\kappa_1(t)} \Pi_n
A_s(u^{n,m}(\kappa_2(s)))\, ds\,  ,\,
\varphi(t)\Big\rangle \, \lambda(t)\, dt                         \nonumber\\
&& \qquad \quad + \sum_{j=1}^r  E\int_0^T\Big(\int_0^{\kappa_1(t)}
 \Pi_n \tilde{B}^{m,j}_s(u^{n,m}(\kappa_1(s)))\, dW^j_s\, ,\,
\varphi(t)\Big) \, \lambda(t) \,dt                                \nonumber\\
&&\quad  =  E\int_0^T(u_0\, ,\,  \varphi(t)) \,  \lambda(t) \, dt +J_1+J_2
-R_1-R_2-R_3  \, ,                                                  \label{23.31.11}
\end{eqnarray}
with
\begin{eqnarray*}
 J_1&:=&E\int_0^T \Big\langle \int_0^t
A_s\big(u^{n,m}(\kappa_2(s))\big)\,ds \, , \,
 \varphi(t)\Big\rangle  \,\lambda(t)\,   dt \, ,\\
 J_2&:=&\sum_{j=1}^r  E\int_0^T\Big(\int_0^t \Pi_{n}
\tilde{B}^{m,j}_s  \big(u^{n,m} (\kappa_1(s))\big) \, dW^j_s\,
,\, \varphi(t)\Big) \,\lambda(t)\,  dt\, ,\\
R_1&:=&E\int_0^{t_1}(u_0\, ,\,  \varphi(t))\,\lambda(t)\,dt\, ,\\
R_2&:=&E \int_0^T\Big\langle \int_{\kappa_1(t)}^t  A_s\big(u^{n,m}(\kappa_2(s))\big)\, ds \, ,\,
 \varphi(t)\Big\rangle  \, \lambda(t)\,  dt \, ,\\
R_3&:=&\sum_{j=1}^r E \int_0^T \Big( \int_{\kappa_1(t)}^t   \Pi_{n}
\tilde{B}^{m,j}_s\big(u^{n,m}(\kappa_1(s))\big)\, dW^j_s \,  ,\,
 \varphi(t)\Big)  \, \lambda(t)\, dt\, .
\end{eqnarray*}

Clearly, for $n,m\to \infty$:
\begin{eqnarray}
|R_1|&\to&0,                                              \label{24.31.11}\\
|R_2|&\leq& N\, E\int_0^T \lambda(t)
\int_{\kappa_1(t)}^t |A_s(u^{n,m}(\kappa_2(s))|_{V^*}\, ds\, dt\nonumber\\
&\leq& N\, 
\Big\{ E\int_0^T  \int_{0}^T\lambda(t) \lambda(s)^{-\frac{q}{p}}\,
|A_s(u^{n,m}(\kappa_2(s))|_{V^*}^q\, ds\,dt
\Big\}^{\frac{1}{q}}\nonumber \\
&&\times
 \Big\{
\int_0^T \lambda(t)\int_{\kappa_1(t)}^t \lambda(s)\,ds\,dt
\Big\}^{\frac{1}{p}}\to 0                         \label{25.31.11} 
\end{eqnarray}  
by 
H\"older's inequality, Lebesgue's theorem 
on dominated convergence, and 
by virtue of estimate (\ref{estimA}).  By the isometry of
$H$-valued stochastic integrals
\begin{eqnarray}
 |R_3|&\leq &
 N \,  \int_0^T E\Big(\int_{\kappa_1(t)}^t \sum_{j=1}^r
\big|\tilde{B}^{m,j}_s\big(u^{n,m}(\kappa_1(s))\big)\big|_H^2\,
ds\Big)^{1/2}\,dt                                      \nonumber\\
&\leq &
  N \,  \sqrt{\delta_m}\, E\int_0^T \sum_{j=1}^r
\big|\tilde{B}^{m,j}_s\big(u^{n,m}(\kappa_1(s))\big)|_H^2\, ds\to
0 \label{22.31.11}
\end{eqnarray}
as $n,m\to \infty$, by virtue of estimate (\ref{estimB}).   
The arguments 
used to prove (\ref{15.30.11})  
in Proposition \ref{uinftyexplicit} yield as $n,m\to\infty$
\begin{equation}
J_{2} 
\to \sum_{j=1}^r
E\int_0^T \Big\langle 
\int_0^t b^j_{\infty}(s)\,d W^j_s,
 \varphi(t)\Big\rangle \,\lambda(t)\,  dt\, .
\end{equation}
Similarly, for $g\in {\mathcal L}^q_{V^*}(\lambda^{1-q})$, 
let $G(g)_t=\int_0^t g(s)\, ds$.
Then H\"older's inequality implies that
\begin{eqnarray*}
\|G(g)\|_{{\mathcal L}^q_{V^*}(\lambda)}^q &\leq & E\int_0^T 
\lambda(t)\, \Big( \int_0^t |g(s)|_{V^*}^q \, \lambda(s)^{-\frac{q}{p}}
\, ds\Big) \, \Big( \int_0^t\lambda(s)\, ds\Big)^{\frac{q}{p}}\, dt \\
&\leq & \Big( \int_0^T \lambda(t)\, dt \Big)^q \, 
\|g\|^q_{{\mathcal L}^q_{V^*}(\lambda^{1-q})}\; .
\end{eqnarray*}
Hence, the operator $G$ is bounded from 
${\mathcal L}^q_{V^*}(\lambda^{1-q})$ to ${\mathcal L}^q_{V^*}(\lambda)$.
Thus this operator is weakly continuous. 
Therefore as $m,n\to \infty$
\begin{equation}
J_1
\to 
E\int_0^T\Big\langle\int_0^t  a_{\infty}(s)
\, d s\, ,\, \varphi(t)\Big\rangle \,\lambda(t)\,  dt\,
.\label{26.31.11}
\end{equation}   
Letting now
$n,m\to\infty$ in equation (\ref{23.31.11}), we obtain
\begin{align*} 
E\int_0^T(v_{\infty}(t)\, ,\, \varphi(t)) &\,
\lambda(t) \, dt =\; E\int_0^T \big(u_0, \varphi(t)\big)\,
\lambda(t)  \, dt 
+ E\int_0^T \Big\langle \int_0^{t}
a_{\infty}(s)\, ds\,  ,\,
\varphi(t)\Big\rangle \, \lambda(t)\, dt      \\*                   
&  + E\int_0^T\Big(\sum_{j=1}^r \int_0^t
b_{\infty}(s)^j\,   dW^j_s\,  ,\,
\varphi(t)\Big) \, \lambda(t) \,dt
\end{align*}
by (\ref{24.31.11})-(\ref{26.31.11}) for any $V_N$-valued adapted
stochastic process $\varphi$ with $\sup_{t,\omega}
|\varphi(t,\omega)|_H  \leq N$.   Since $N$ can be arbitrary large,
equation (\ref{evolimplicit}) follows immediately.
As in the proof of (\ref{16.29.11}), a similar argument based on
an analog of (\ref{22.30.11}) for a $L^2_{V_N}$ random variable
$\psi$ with $E|\psi|_V^2\leq N$ yields  equation (\ref{16.31.11}).
An argument similar to that proving (\ref{22.31.11}) yields 
\begin{equation}                                      \label{22.01.11}
E|u_{\infty}(T)|_H^2 =E|u_0|^2_H+E\int_0^T  \Big[ 2\langle
v_{\infty}(s),a_{\infty}(s)\rangle +\sum_{j=1}^r
|b_{\infty}(s)|_H^2 \Big] \,ds\; . 
\end{equation}
Moreover, by (\ref{evolimplicit}) and (\ref{16.31.11}) we get 
$u_{\infty}(T)=u_{\infty T}$ (a.s.).  To prove inequality (\ref{20.01.11}) set
\begin{eqnarray*}
&&
F^{n,m}_y:=E\int_0^T2\Big\{\, \big\langle u^{n,m}(\kappa_2(t))-y(t)\, ,
\,  A_t(u^{n,m}(\kappa_2(t))) -
A_t(y_t)\big\rangle
\\ &&\qquad \qquad \quad +\sum_{j=1}^r  \big|B_t^j(u^{n,m}(\kappa_2(t)))
-\Pi_nB_t^j(y(t))\big|_H^2\Big\}\,dt
\end{eqnarray*}
for $y\in{\mathcal L}^p_{V} (\lambda)\cap {\mathcal L}^2_H(K_1)$.  By (C4),
(\ref{10.30.11}) and Lemma \ref{apriori}$, F^{n,m}_y$ is
well-defined and it is finite.  By the monotonicity condition and
by inequality (\ref{equal}) with $k:=m$  we obtain:
  \begin{align}
0\geq  F^{n,m}_y& 
\geq\;  E|u^{n,m}(T)|_H^2 - E|u_0|_H^2 +2\,  E\int_0^T\langle y_t
\, \, A_t(y_t)\rangle \, dt   -2\,L_1^{n,m} \nonumber\\*
  & - 2\, L_2^{n,m}+L_3^{n}-2\, L_4^{n,m}
 + \delta_m\, E \int_{\delta_m}^T \big|
\Pi_{n} A_s(u^{n,m}(\kappa_2(s)))\, \big|_H^2\, ds\, ,
\label{ineqyk}
\end{align}
with
\begin{eqnarray*}
L_1^{n,m}&:= & E\int_0^T \big\langle u^{n,m}(\kappa_2(t))\,
, A_t(y_t)\big\rangle \, dt      \, ,    \\ 
L_2^{n,m}&:=&E\int_0^T \big\langle y_t\, , \,
A_t\big(u^{n,m}(\kappa_2(t))\big)\big\rangle \, dt\, ,        \\ 
L_3^{n}&:=&\sum_{j=1}^rE\int_0^T \big|\Pi_{n}B_t^j(y_t)\big|_H^2\,dt \, ,\\
L_4^{n,m}&:=& \sum_{j=1}^r
E\int_0^T\Big(\Pi_{n}B_t^j\big(u^{n,m}(\kappa_2(t))\big)\, ,\,
B_t^j(y_t)\Big)\,  dt\, .
\end{eqnarray*}
Using (i)-(iv) and the arguments used to prove (\ref{19.30.11}),
(\ref{20.30.11}),  (\ref{17.30.11bis}) and (\ref{21.30.11}) we deduce:
\begin{eqnarray}
\lim_{n,m\to\infty}L_1^{n,m}&=&E\int_0^T\langle  v_\infty(t)  \, ,\,
A_t(y_t)\rangle\, dt,  \label{11.03.11}\\
\lim_{n,m\to\infty}L_2^{n,m}&=&E\int_0^T\langle y_t\, , \,
a_{\infty}(t)\rangle\,dt, \label{10.03.11}\\
\lim_{n\to\infty}L_3^{n}& =& \sum_{j=1}^r \int_0^T|B_t^j(y_t)|_H^2 \,
dt\,, \\
\lim_{n,m\to\infty}L_4^{n,m}&=&\sum_{j=1}^r E\int_0^T
\big(b_{\infty}^j(t)\, ,\, B^j_t(y_t)\big)\,dt. \label{13.03.11}
\end{eqnarray}
Furthermore, for some constant $d\geq 0$:
\begin{equation}
\liminf_{n,m\to\infty}E|u^{n,m}(T)|_H^2 =d+E|u_{\infty}(T)|_H^2.
\label{09.03.11}
\end{equation}
Thus, letting $n,m\to\infty$ in  (\ref{ineqyk}), by
(\ref{11.03.11})-(\ref{09.03.11}) 
we deduce:
\begin{eqnarray*}
 0&\geq& d + E|u_\infty(T)|_H^2  - E|u_0|_H^2 -2E\int_0^T
\langle  v_\infty(t)\, ,\, A_t(y_t)\rangle\, dt 
-2 E\int_0^T\langle y_t\, ,\, a_\infty(t) \rangle\, dt 
\\&&
 + 2\, E\int_0^T\langle y_t\, ,\, A_t(y_t)\rangle\, dt
 + \sum_{j=1}^r  E\int_0^T \big[ |B_t^j(y_t)|_H^2
-2\, (b_\infty^j(t), B_t^j(y_t))\big]\, dt                        \\
& =&d+F_y(u_0,a_{\infty},b_{\infty})\, .
\end{eqnarray*}
Then we proceed as after (\ref{inegalex}) at the end
of the proof of Proposition \ref{exp} and finish
the proof of the proposition.
\end{proof}
\bigskip

Now we conclude the proof of Theorem \ref{convergence}.   By the
previous proposition, from any sequence $(n,m)$ of pairs of
positive integers such that $m,n \to \infty$, there exists a
subsequence, $(n_k,m_k)$, such that the approximations
$u^{n_k,m_k}$ converge weakly in ${\mathcal L}^p_V(\lambda)$ to
the solution $u$ of equation (\ref{u}), and the approximations
$u^{n_k,m_k}(T)$ converge strongly in $L^2_H$ to $u(T)$.   Hence,
taking into account that the solution of equation (\ref{u}) is
unique, we get that these convergence statements hold for any
sequences of approximations $u^{n,m}$ and $u^{n,m}(T)$ as
$n,m \to\infty$.   The proof of Theorem \ref{convergence} is
complete.    \hfill $\Box$

\bigskip

\noindent {\bf Acknowledgments:}
The authors wish to thank the referee for
helpful comments.

\section{Appendix}

We start with a technical lemma ensuring that a map
from $V$ to $V^*$ is continuous.
\begin{lem}\label{continuity}
Let $V$ be a Banach space and $V^*$ its topological dual,
$D:V\rightarrow V^*$  satisfy the conditions (i)-(iii) of
Proposition \ref{infdim}.  Then $D$ is continuous from $(V,|\, .\,
|_V)$ into $V^*$ endowed with the weak star topology
$\sigma(V^*,V)$.   In particular, if $V$ is a  finite dimension
vector space, then $D$ is continuous.
\end{lem}
\begin{proof} Let $x\in V$ and $(x_n,\,  n\geq 1)$ be a sequence of
elements of $V$ such that $\lim_n|x-x_n|_V=0$.
The monotonicity property (i) implies that for every $y\in V$ and $n\geq 1$,
$$\langle D(x_n)-D(y), x_n-x\rangle +\langle D(x_n)-D(y), x-y\rangle  \geq 0\, .
$$
Furthermore, since  $(|x_n|_V\, ,\, n\geq 1)$ is bounded,
the growth condition (iii) implies that
\begin{eqnarray*}
|\langle D(x_n)-D(y),x_n-x\rangle |
& \leq& \big[ |D(x_n)|_{V^*}+|D(y)|_{V^*}\big]\, |x_n-x|_V \\
& \leq & C \,(1+|x_n|_V^{p-1}+|y|_V^{p-1})\,  |x_n-x|_V\rightarrow 0
\end{eqnarray*}
as $n\rightarrow +\infty$; hence,
$\liminf_n \langle D(x_n)-D(y),x-y\rangle  \geq 0$.
Since $(|x_n|_V, n\geq 1)$
is bounded, the growth condition implies the existence of
a subsequence $(n_k , k\geq 1)$ such that
$D(x_{n_k}) \rightarrow D_\infty\in V^*$ is the weak
star topology as $k\rightarrow +\infty$; clearly,
\begin{equation}\label{ineq}
\langle D_\infty -D(y),x-y\rangle  \geq 0\; , \forall  y\in V\, .
\end{equation}
To conclude the proof, we check that $D_\infty = D(x)$; indeed,
this yields that the whole sequence $(D(x_n), n \geq 1)$ converges
weakly to $D(x)$.   For any $z\in V$ and $\varepsilon >0$, apply
(\ref{ineq}) with $y=x - \varepsilon z$; then dividing by
$\varepsilon >0$ and using the hemicontinuity property (ii), we
deduce that for any $z\in V$,
$$
\lim_{\varepsilon \rightarrow 0}\langle D_\infty - D(x-\varepsilon
z),z\rangle = \langle D_\infty -D(x),z\rangle \geq 0\, .
$$
Changing $z$ into $-z$, this  yields $D_\infty = D(x)$.
\end{proof}

{\it Proof of Proposition \ref{infdim}.}
Let $(e_i\, ,\, i\geq 1)$ be a sequence of elements  of $V$ which is 
a complete orthonormal basis of $H$ and for every $n\geq 1$, 
let $\tilde{V}_n= \mbox{\rm span }(e_i\, ,\, 1\leq i\leq n)$,
 $\tilde{\Pi}_n:V^*\rightarrow \tilde{V}_n$ 
be defined by $\tilde{\Pi}_n(y)=\sum_{i=1}^n \langle e_i,y\rangle \, e_i$ for
 $y\in V^*$ and let $\tilde{D}_n= \tilde{\Pi}_n\circ D :\tilde{V}_n \rightarrow \tilde{V}_n$. 
Then $\tilde{D}_n$ is coercive and satisfies the assumptions of Lemma \ref{continuity}; 
hence it is continuous.   Fix $y\in V^*$;
 the existence of $x_n\in \tilde{V}_n$ such that 
$\tilde{D}_n(x_n)=\tilde{\Pi}_n(y)$ is classical (see e.g. \cite{Z}).
The coercivity condition implies that for every $n\geq 1$:
 \[ |y|_{V^*}\, |x_n|_V\geq \langle x_n,y\rangle 
= \langle x_n,\tilde{D}_n(x_n)\rangle =  \langle x_n, D(x_n)\rangle
 \geq C_1\, |x_n|^p_{V} - C_2\, ,\]
which implies that the sequence $(|x_n|_V, n\geq 1)$ is bounded, 
and the growth property implies that the sequence
 $(|D(x_n)|_{V^*}\, ,\, n\geq 1)$ is boun\-ded. 
Since $V$ is reflexive, there exists a subsequence  $(n_k\, ,\, k\geq 1)$ such
 that the sequence  $(x_{n_k}\, ,\, k\geq 1)$ converges to  $x_\infty \in V$ 
in the weak $\sigma(V,V^*)$ topology, and such that the sequence 
$(D(x_{n_k})\, ,\, k\geq 1)$ converges to $D_\infty$ in the weak-star topology
 $\sigma(V^*, V)$.   We at first check that $D_\infty = y$; indeed, for every $i\geq 1$:
 \begin{eqnarray*} \langle D_\infty, e_i\rangle &=&\lim_k \langle D(x_{n_k}), e_i\rangle=
 \lim_k \langle \tilde{D}_{n_k}(x_{n_k}), e_i\rangle \\
&=& \lim_k \langle \tilde{\Pi}_{n_k}(y),e_i\rangle
 = \langle y,e_i\rangle\, . 
\end{eqnarray*}
We then prove that $y=D(x_\infty)$; the monotonicity property of $D$ implies 
that for every $z\in \cup_n V_n $, for $k$ large enough:
 \begin{eqnarray*}
 0&\leq &\langle D(x_{n_k})-D(z),x_{n_k} -z\rangle\\
 &\leq & \langle \tilde{D}_{n_k}(x_{n_k}), x_{n_k}\rangle- \langle D(z), x_{n_k}\rangle
 - \langle \tilde{D}_{n_k}(x_{n_k}),z \rangle  + \langle D(z),z\rangle\\
 &\leq& \langle y, x_{n_k}\rangle - \langle D(z), x_{n_k}\rangle 
- \langle \tilde{D}_{n_k}(x_{n_k}),z \rangle
 + \langle D(z),z\rangle\, .
 \end{eqnarray*}
As $k\rightarrow +\infty$, we deduce that for every 
$z\in \cup_n V_n$, $\langle D_\infty-D(z), x_\infty-z\rangle \geq 0$. 
Since $\cup_n V_n$ is dense in $V$, we deduce that 
$\langle D_\infty-D(z), x_\infty-z\rangle \geq 0$ for every $z\in V$. 
Let $\xi \in V$;  apply the previous inequality to $z=x_\infty +\varepsilon\, \xi$  
for any $\varepsilon >0$ and divide by $\varepsilon$. 
This yields that for any $\xi \in V$, $\langle D_\infty - D(x_\infty +
\varepsilon \, \xi)\, ,\,\xi\rangle \geq 0$; as $\varepsilon \rightarrow 0$, 
the hemicontinuity implies that for any $\xi\in V$,
$\langle D_\infty - D(x_\infty)\, ,\,  \xi
\rangle\geq 0$, 
and hence that $y=D_\infty = D(x_\infty)$.   This concludes the proof of the 
existence of a solution $x$ to the equation $D(x)=y$. 
Furthermore, the coercivity of $D$ implies that
\[C_1\, |x|_V^p -C_2 \leq \langle D(x)\, ,\, x \rangle 
= \langle y\, ,\, x\rangle
 \leq \frac{C_1}{2}\, |x|_V^2 + \frac{1}{2 C_1}\, |y|_{V^*}^2\, . \] 
Hence for $p\in [2,+\infty[$, $C_1\, |x|_V^p -C_2 \leq
\frac{C_1}{2}\, |x|_V^p + \frac{C_1}{2}+ \frac{1}{2 C_1}\, |y|_{V^*}^2$,
which implies (\ref{norm}).   Finally, if  $D$ satisfies the strong monotonicity 
condition (\ref{strongmono}) and if $x_1, x_2\in V$ 
are such that $D(x_1)=D(x_2)=y$, then
\[ 
0=\langle D(x_1)-D(x_2), x_1-x_2\rangle \geq C_3\, |x_1-x_2|^2_{V^*}\, ;
\]
this  yields $|x_1-x_2|_{V^*}=0$.\quad $\Box$
\bigskip

We finally sketch the proof of Theorem \ref{minimization}

{\it Proof of Theorem  \ref{minimization}.} (i) The
monotonicity condition (C1) implies that for every $y\in {\mathcal L}^p_V(\lambda)$
such that ${\mbox {\rm sup}_{0\leq t\leq T} } E|y_t|_H^2 <+\infty$ one has:
\[ 
2\langle u_s-y_s\, ,\, A_s(u_s)-A_s(y_s)\rangle
+\sum_{j=1}^r |B_s^j(u_s)-B_s^j(y_s)|_H^2 \leq 0\, .
\] 
This implies that $F_y(u_0, A_.( u_.), B_.(u_.)) \leq 0$ for every
$y\in {\mathcal L}^p_V(\lambda)$ with $\sup_{0\leq t\leq T}
E|y_t|_H^2<+\infty$, which yields (i).

(ii) Let $(\xi,a,b)\in {\mathcal A}$, $u_t=\xi+\int_0^t a_s\, ds +
\sum_{j=1}^r \int_0^t b_s^j\, dW^j_s$ and let ${\mathcal V}$ be a
  subset of ${\mathcal L}^p_V(\lambda)$ of processes $y$ such that
  $\sup_{0\leq t\leq T} E|y_t|_H^2<+\infty$, 
which is dense in ${\mathcal L}^p_V(\lambda)$ and  such that
\begin{equation}                                            \label{24.03.11}
F_y(\xi,a,b)\leq 0
\quad {\mbox {\rm for }} y\in {\mathcal V}\,.
\end{equation}
We first check that (\ref{24.03.11}) holds for $y=u+z$ where
$\sup_{0\leq t_leq T}E|z_t|_V^p<+\infty$.   To this end let $\{y_n\,
,\,n\geq 1\}$ be a sequence of elements of ${\mathcal V}$, such
that $\lim_n \| y-y_n\|_{{\mathcal L}^p_V(\lambda)} =0$.
 For any $U\in {\mathcal L}^p_V(\lambda)$  such that 
${\mbox{\rm sup}_{0 \leq t \leq T}  }   E|U_t|_H^2<+\infty$, 
set
\[ \Phi(U)=E\int_0^T \Big[ 2\big\langle u_s-U(s)\,
,\, a_s-A_s(U(s))\big\rangle +\sum_{j=1}^r
\big|b_s^j-B_s^j(U(s))\big|_H^2\Big]\, ds\, .\]
Then $\big|\, \Phi(y_n)-\Phi(y)\Big| \leq \sum_{i=1}^3 T_i(n)$, where:
\begin{eqnarray*}
T_1(n)&=&\Big| E\int_0^T 2\, \big\langle y(s)-y_n(s)\,
,\, a_s-A_s(y_n(s))\big\rangle\, ds\Big|\, ,\\
T_2(n)&=& \Big| E\int_0^T 2\, \big\langle u_s-y(s)\,
,\, A_s(y(s))-A_s(y_n(s))\big\rangle\, ds\Big|\, ,\\
T_3(n)&=& \sum_{j=1}^r
\Big| E\int_0^T \Big[ \big|B_s^j(y_n(s))\big|_H^2 - \big|B_s^j(y(s))\big|_H^2
+ 2\, \big( b_s^j\, ,\,
B_s^j(y(s))- B_s^j(y_n(s))\big)\Big]\, ds \Big|\, .
\end{eqnarray*}
Since $\sup_n E\int_0^T|y_n(s)|_V^p\, \lambda(s) \, ds<\infty$,
the growth condition (C4) yields
 \begin{eqnarray}
T_1(n)&\leq&\|y-y_n\|_{{\mathcal L}_V^p(\lambda)} \,
\Big\{E \int_0^T \big( |a_s|^q_{V^*}
 + |A_s(y_n(s))|_{V^*}^q\big) \lambda^{1-q}(s)\, ds \Big\}^{\frac{1}{q}}          \nonumber\\
&\leq& C_1\, \|y-y_n\|_{{\mathcal L}_V^p(\lambda)} 
\Big\{ E\int_0^T \Big[ |a_s|_{V^*}^q\, \lambda^{1-q}(s)
 +\alpha\,   |y_n(s)|_{V}^p\, \lambda(s)  +
K_2(s)\Big]\, ds  \Big\}^{\frac{1}{q}}\nonumber \\
&\leq& C_2 \, \|y-y_n\|_{{\mathcal L}^p_V(\lambda)} \, ,              \label{T1n}
\end{eqnarray}
where $C_1,C_2$ are constants which do not depend on $n$.
For $dt\times P$-almost every $(t,\omega)$ the operator
$A_t(\omega)\, :\, V\rightarrow V^*$ is monotone and
hemicontinuous, hence it is demi-continuous, i.e.,
the sequence $A_t(\omega,x_n)$ converges weakly in $V^*$ to $A_t(\omega,x)$
whenever $x_n$ converges strongly in $V$ to $x$ (see, e.g., Proposition 26.4 in \cite{Z} ).
Hence for $ dt\times P-$almost every $(t,\omega)\in[0,T]\times\Omega$, 
$$
\lim_n \langle z(s)\, ,\, A_s(y(s))- A_s(y_n(s))\rangle =0\; .
$$
Furthermore, since $z$ is bounded, condition (C4) implies
\begin{align*}
 \sup_n E\int_0^T \big|\big\langle & z(s)\, ,
\, A_s(y(s))-A_s(y_n(s))\big\rangle \big|^q\, ds\\
&\;\; \leq C \sup_n E\int_0^T\!|A_s(y(s))- A_s(y_n(s))|_{V^*}^q \, ds \\
&\; \; \leq C_1 \sup_n E\int_0^T\! \big(|y(s)|_V^p + |y_n(s)|_V^p\big)\, \lambda(s) \, ds
+ C_1\, \int_0^TK_2(s)\, ds <\infty \, .
\end{align*}
Therefore, the sequence
$\big\{\big\langle z\, ,\, A(y)-A(y_n)\big\rangle\, ,\, n\geq 1 \big\}$
is uniformly integrable with respect to the measure $dt\times P$.   Hence
\begin{equation}\label{T2n}
\lim_n T_2(n)=0\, .
\end{equation}
By Remark \ref{15.29.10}
\begin{eqnarray*}
&&\sup_n\sum_{j=1}^r  E\int_0^T  \big[ |B_s^j(y(s))|_H^2 + |B_s^j(y_n(s))|_H^2 \big]\, ds\\
&&\qquad\quad \leq \sup_n \, CE\int_0^T \big[|y(s)|_V^p + |y_n(s)|_V^p\big]\,\lambda(s)\, ds  \\
&&\quad \qquad\qquad + C\, E\int_0^T \big\{K_1(s)\,
\big[|y(s)|^2_H+|y_n(s)|_H^2\big] +K_3(s)\big\}\, ds <\infty\, .
\end{eqnarray*}
By Schwarz's inequality   we deduce
\begin{eqnarray*}
T_3(n)&\leq & C\, \sum_{j=1}^r \Big\{ E\int_0^T  \big[
|B_s^j(y(s))|_H^2+ |B_s^j(y_n(s))|_H^2+  |b^j_\infty(s)|_H^2\big]
\, ds \Big\}^{\frac{1}{2}}\\
&&\qquad \times \Big\{ E\int_0^T
|B_s^j(y_n(s)) - B_s^j(y(s))|_H^2\, ds \Big\}^{\frac{1}{2}}\\
& \leq& C\, \sum_{j=1}^r \Big\{ E\int_0^T|B_s^j(y_n(s)) - B_s^j(y(s))|_H^2\,
ds \Big\}^{\frac{1}{2}}\, .
\end{eqnarray*}
The monotonicity assumption (C1) and the growth condition (C4) imply:
\begin{align}
\nonumber \sum_{j=1}^r&  E\int_0^T\! \!
\big|B_s^j(y_n(s)) - B_s^j(y(s))\big|_H^2\, ds
\leq  -2\, E\int_0^T\!\!
\big\langle y_n(s)-y(s)\, ,\,
A_s(y_n(s))-A_s(y(s))\big\rangle\, ds   \nonumber \\
&\leq C\, \Big\{ E\int_0^T\!\!  |y_n(s)-y(s)|_V^p\,
\lambda(s)\,  ds \Big\}^{\frac{1}{p}} 
\Big\{E \int_0^T \!\! \big[ \big( |y_n(s)|_V^p+|y(s)|_V^p\big )\, \lambda(s)
 + K_2(s)\big] ds  \Big\}^{\frac{1}{q}}  \nonumber \\
&\leq C\, \| y_n-y\|_{{\mathcal L}_V^p(\lambda)}\, . \label{T3n}
\end{align}
The inequalities (\ref{T1n})-(\ref{T3n}) imply $\lim_n \Phi(y_n)=\Phi(u+z)$.  
Consequently, (\ref{24.03.11}) holds for $y=u+z$ with any
$z\in{\mathcal L}^{\infty}_V$.

Fix $z\in {\mathcal L}^{\infty}_V$ and $\varepsilon >0$, apply
(\ref{24.03.11}) to $y=u -\varepsilon z$ and divide by
$\varepsilon$; this yields:
\begin{equation}                                        \label{23.03.11}
 E \int_0^T   \big\langle z_t, a_t-A_t(u_t-\varepsilon z_t)
\big\rangle \, dt \geq 0\, .
\end{equation}
By the hemicontinuity condition (C3) for almost all $\omega\in\Omega$:
$$
\lim_{\varepsilon\to0}\big\langle z_t\, ,\, a_t-A_t(u_t-\varepsilon z_t)\big\rangle
=\big\langle z_t\, ,\, a_t-A_t(u_t)\big\rangle\, , \quad\forall t\in[0,T].
$$
Furthermore, since $z$ is bounded, by (C4)
\[
\sup_{0<\varepsilon \leq 1} E\int_0^T
|\langle z_t\, ,\, a_t-A_t(u_t-\varepsilon\, z_t)\rangle |^q\, dt
<\infty\, ,
\]
which implies that
$\{\langle z\, ,\, a-A(u-\varepsilon z)\rangle\, ,\, 0<\varepsilon \leq1\}$
is uniformly integrable over $[0,T]\times\Omega$, with respect to the measure 
$dt\times P$.  Hence letting $\varepsilon \rightarrow 0$ in (\ref{23.03.11}) we get
\[
 E\int_0^T \left\langle  z_t\, ,\, a_t-A_t(u_t)) \right\rangle\, dt
\leq 0 \quad \mbox{\rm for any}
 \quad z\in {\mathcal L}^{\infty}_V \, .
\]
Changing $z$ into $-z$ and using that ${\mathcal L}^{\infty}_V$ is dense in
${\mathcal L}^{p}_V(\lambda)$ we deduce that
 \[
a_t(\omega)=A_t(u_t(\omega),\omega)\quad {\mbox{\rm for }}
dt\times P{\mbox{\rm almost every }}
(t,\omega)\in[0,T]\times\Omega.
\]
Using again (\ref{24.03.11}) with $y=u$ (i.e., $z=0$), we deduce
that $B_t(u_t(\omega),\omega)$ $=b_t(\omega)$ for $dt\times P$
almost every $(t,\omega)$, and that $\xi=u_0$ (a.s.).
Consequently, $u$ is a solution to (\ref{u}). \qquad $\Box$

\end{document}